\input amstex
\input epsf
\documentstyle{amsppt}
\NoBlackBoxes
\magnification 1200
\vcorrection{-13mm}  

\topmatter
\title     Markov trace on the Funar algebra
\endtitle

\author    S.~Yu.~Orevkov
\endauthor

\address
           IMT, Universit\'e Paul Sabatier, Toulouse, France
\endaddress
\address
           Steklov Mathematical Institute, Moscow, Russia
\endaddress

\endtopmatter

\document

\def\refBF               {1}
\def\refCM               {2}
\def\refFunar            {3}
\def\refMarin            {4}
\def\refO                {5}
\def\refWebPage          {6}
\def\refOS               {7}

\def\eqRelOne            {1}
\def\eqRelTwo            {2}
\def\eqRelTwoBar         {3}
\def\eqMT                {4}
\def\eqRelGeneral        {5}
\def\eqDefS              {6}
\def\eqTau               {7}
\def\eqWt                {8}
\def\eqY                 {9}
\def\eqNonInteractOne   {10}
\def\eqNonInteractTwo   {11}
\def\eqInterB           {12}
\def\eqSkein            {13}
\def\eqDefP             {14}

\def\propReduc    {2.1}
\def\remOne       {2.2}
\def\remTwo       {2.3}
\def\thMain       {2.4}
\def\corZZZ       {2.5}
\def\corBetaZero  {2.6}

\def\remInv       {2.8}
\def\remMarin     {2.9}
\def\remTransv    {2.10}
\def\remFunarErr  {2.11}

\def\lemEasy       {3.1}
\def\propEasy      {3.2}
\def\lemDiffiOne   {3.3}
\def\lemDiffiTwo   {3.4}
\def\lemPent       {3.5}
\def\lemDiffiThree {3.6}
\def\propDiffi     {3.7}

\def\propInv       {4.1}
\def\propPeriod    {4.2}

\def\propGBQ       {4.3}
\def\propGBZ       {4.4}
\def\remMonomOrder {4.5}
\def\propSym       {4.6}
\def\examFigEight  {4.7}
\def\propHOMFLY    {4.8}
\def\examHOMFLY    {4.9}
\def\propSpin      {4.10}

\def\sectReduc       {2.1}
\def\sectMarkov      {2.2}
\def\sectMainResult  {2.3}

\def\sectInv    {4.1}
\def\sectNF     {4.2}
\def\sectGBQ    {4.3}
\def\sectGBZ    {4.4}
\def\sectHOMFLY {4.5}

\def\figBunch       {1}
\def\figNonInteract {2}
\def\figInteract    {3}
\def\figInteractA   {4}
\def\figInteractB   {5}

\def\Z{\Bbb Z}
\def\eps{\varepsilon}

\def\red{\operatorname{red}}

\def\sh{\operatorname{sh}}

\def\wt{\operatorname{wt}}

\def\pr{\operatorname{pr}}
\def\End{\operatorname{End}}
\def\Hom{\operatorname{Hom}}
\def\rr{\bold{r}}
\def\tto{\longrightarrow}
\def\u{\bar u}
\def\v{\bar v}


\hskip 72mm\epsfxsize=47mm\epsfbox{epi2500.eps}
\medskip

\head 1. Introduction
\endhead

Let $B_n$ be the braid group with $n$ strings and 
$\sigma_1,\dots,\sigma_{n-1}$
its standard generators.
Let $k$ be a commutative ring with $1\ne 0$.
Given $\alpha,\beta\in k$, we define the $k$-algebra
$K_n=K_n(\alpha,\beta)=K_n(\alpha,\beta;k)$ as the quotient
of the group algebra $kB_n$ by the relations
$$
     \sigma_1^3-\alpha\sigma_1^2+\beta\sigma_1-1=0  \eqno(\eqRelOne)
$$
and
$$
\split
  y\bar xy=&2\alpha-\beta^2-(x+y)-(\alpha^2-\beta)(\bar x+\bar y)
  +\beta(xy+yx) + \alpha(x\bar y+y\bar x+\bar xy+\bar yx)
\\
  &+(\alpha\beta-1)(\bar x\bar y+\bar y\bar x)
 - \alpha xyx
  -(\bar xyx+x\bar yx+xy\bar x) - \beta(\bar x\bar yx +x\bar y\bar x)
\\
   &+ (\alpha-\beta^2)\bar x\bar y\bar x.
  \hskip 88mm                                            (\eqRelTwo)
\endsplit
$$
where $x$, $\bar x$, $y$, $\bar y$ in (\eqRelTwo) stand for
$\sigma_1$, $\sigma_1^{-1}$, $\sigma_2$, $\sigma_2^{-1}$ respectively.
Up to a change of the sign of $\beta$ (for the sake of symmetricity), our
definition of $K_n$ is equivalent to the definition given by
Bellingeri and Funar in [\refBF]. Our relation (\eqRelTwo) is much shorter
than the corresponding relation in [\refBF] (see [\refBF; (2) and Table 1])
because we use $\sigma_i^{-1}$ instead of
$\sigma_i^2$. Multiplying (\eqRelTwo) by $\sigma_1$ from the left or
from the right, and simplifying the result using (\eqRelOne) and the braid group
relations, we obtain 
$$
\split
  \bar y x\bar y=&2\beta-\alpha^2-(\bar x+\bar y)-(\beta^2-\alpha)(x+y)
  +\alpha(\bar x\bar y+\bar y\bar x) + \beta(\bar xy+\bar yx+x\bar y+y\bar x)
\\
  &+(\alpha\beta-1)(xy+yx)
 - \beta \bar x\bar y\bar x
  -(x\bar y\bar x+\bar x y\bar x+\bar x\bar y x) - \alpha(xy\bar x +\bar xyx)
\\
   &+ (\beta-\alpha^2)xyx
  \hskip 88mm                                            (\eqRelTwoBar)
\endsplit
$$
(note that (\eqRelTwoBar) is obtained from (\eqRelTwo) by swapping
$x\leftrightarrow\bar x$,
$y\leftrightarrow\bar y$,
$\alpha\leftrightarrow\beta$).

%
Using (\eqRelOne) -- (\eqRelTwoBar) together with the braid
relations, it is easy to see that $K_n$ are finitely generated $k$-modules.
Following [\refCM], we denote the image of $\sigma_i$ in $K_n$ by $s_i$.

Set $K_\infty=\lim K_n$ (in contrary to the case of Hecke or BMW algebras,
the morphisms
$K_n\to K_{n+1}$ induced by the standard embeddings $B_n\subset B_{n+1}$
are not injective in general).
We say that $t:K_\infty\otimes k[u,v]\to M$ is a {\it Markov trace} on $K_\infty$ if
$M$ is a $k[u,v]$-module and $t$ is a morphism of $k[u,v]$-modules such that 
$t(xy)=t(yx)$,
$t(xs_n)=ut(x)$,
$t(xs_n^{-1})=vt(x)$, $x,y\in K_n$, $n=1,2,\dots$.

It is claimed in [\refFunar] and [\refBF] that a nontrivial
Markov trace is constructed on $K_n$.
About 2004--2005 I indicated a gap in the proof of its well-definedness
(see Remark \remFunarErr\ below).
As it is explained in [\refCM],
the gap was really serious: formally, the main result of [\refFunar]
is wrong in the form it is stated.
However, we show in this paper that
the main idea in [\refBF, \refFunar]
is correct: to construct a Markov trace on $K_n$, it suffices to check a finite
number of identities though the number of them is much
bigger than in [\refBF, \refFunar] and the algorithm of computation
is much more complicated.
Theoretically, this approach allows to compute the universal Markov trace on $K_\infty$,
i.~e., the projection of $K_\infty(\alpha,\beta;\Z[\alpha,\beta,u,v])$ onto
its quotient by the submodule $\bar R$
generated by 
$$
   xy-yx,\;\;
   xs_n-ux,\;\;
   xs_n^{-1}-vx,\quad x,y\in K_n,\; n=1,2,\dots
                                                                     \eqno(\eqMT)
$$
The volume of computations is
huge, so we performed them only in some cases.
In particular, we found 
$K_\infty(0,0;A)/\bar R=A/I$
where $A=\Z[u,v]$, $I=(16,4u^2+4v,4v^2+4u,u^3+v^3+uv-3)$.
Note, that it was checked in [\refCM] that 
$K_5(0,0;A)/\bar R_5=A/I$ where $\bar R_5$ is the submodule
generated by the elements of $K_5$ of the form (\eqMT).

In a sense, the results of the present paper can be divided into two independent parts:
the ``theoretical part'' (Theorem \thMain\ which provides an algorithm for computing the ideal $I$)
and the ``computational part'' (Corollaries \corZZZ\ and \corBetaZero\ which present
the results of computer-aided computations according to
this algorithm, and Section 4 where we discuss some properties of the obtained link invariants). 
The explicit form of the coefficients in the right hand side of (\eqRelTwo) and (\eqRelTwoBar) is not
really used in the
``theoretical part''. Theorem \thMain\ can be applied to a quotient of $kB_\infty$
by (\eqRelOne) together with any two relations of the form
$$
   \sigma_2\sigma_1^{-1}\sigma_2 = \sum_{i=1}^{21} \gamma_i X_i,
   \qquad
   \sigma_2^{-1}\sigma_1\sigma_2^{-1} = \sum_{i=1}^{21}
            \gamma'_i X_i,\qquad \gamma_i,\gamma'_i\in k,                        \eqno(\eqRelGeneral)
$$
where $X_1,\dots,X_{21}$ are all the reduced words in $\sigma_1^{\pm1}$, $\sigma_2^{\pm1}$
(including the empty word $1$)
that do not contain any subword of the form $\sigma_i^{\pm2}$ or
$\sigma_2^{\pm1}\sigma_1^{\pm1}\sigma_2^{\pm1}$.
So, if we consider $\alpha,\beta$ and all the $\gamma_i$, $\gamma'_i$ in
(\eqRelGeneral) as
indeterminates and compute the ideal $I$ described in Theorem \thMain,
then we obtain the universal Markov trace on a cubic Hecke algebra that can be specialized
to both Funar and BMW algebras.
However, the required computations seem to exceed the capacity of any computer.
On the other hand, if one chooses a random specialization of the coefficients
$\gamma_i$, $\gamma'_i$ in (\eqRelGeneral),
then the quotient $A/I$ might well be trivial.
So, the explicit form of (\eqRelTwo) and (\eqRelTwoBar) is
important for the ``computational part'' of the present paper.


\subhead Acknowledgements
\endsubhead
%
I am grateful to Andrey Levin and Alexey Muranov for useful discussions and advice.
I am grateful to the referees for valuable remarks.

\head 2. Definitions and statement of results
\endhead

\subhead\sectReduc. $K$-reductions
\endsubhead
Let $F_n^+$ be the free monoid on generators
$x_1^{\pm1},\dots,x_{n-1}^{\pm1}$ (the set
of all not necessarily reduced words in $x_i^{\pm1}$)
and $F_\infty^+=\bigcup F_n^+$. We denote the empty word by $1$.
Let $kF_n^+$ and $kF_\infty^+$ be the corresponding free associative
algebras over $k$ (as $k$-modules, they are freely generated by $F_n^+$
and by $F_\infty^+$ respectively).

We call {\it basic replacements} the pairs $(U,V)$ with $U\in F_\infty^+$,
$V\in kF_\infty^+$ (which we denote by $U\to V$) from the following list:
   \smallskip

\roster
\item "$(i)$"
   $x_i x_i^{-1}\tto 1$, $\;x_i^{-1} x_i\tto 1$, $\,i\ge1$;
   \smallskip
\item "$(ii)$"
   $x_i^2\tto \alpha x_i - \beta + x_i^{-1}$, $\,i\ge 1$;
   \smallskip
\item "$(iii)$"
   $x_i^{-2}\tto \beta x_i^{-1} - \alpha +  x_i$, $\,i\ge 1$;
   \smallskip
\item "$(iv)$"
   $ x_{i+1}^{\eps_1}\,  x_i^{\eps_2}\,  x_{i+1}^{\eps_3} \tto
     x_i^{\eps_3}  x_{i+1}^{\eps_2}  x_i^{\eps_1}$,
          $\,\eps_2\in\{\eps_1,\eps_3\}\subset\{-1,1\}$, $\,i\ge 1$;
   \smallskip
\item "$(v)$"
   $ x_{i+1} x_i^{-1} x_{i+1} \tto$
   (the right hand side of (\eqRelTwo) with $x= x_i$, $y= x_{i+1}$), $\,i\ge1$;
   \smallskip
\item "$(vi)$"
   $\, x_{i+1}^{-1}\, x_i\, x_{i+1}^{-1}\, \tto$
   (the right hand side of (\eqRelTwoBar) with $x= x_i$, $y= x_{i+1}$), $\,i\ge1$;
   \smallskip
\item "$(vii)$"
   $ x_{i+1}^{\eps_1} x_i^{\eps_2} W  x_{i+1}^{\eps_3} \tto VW$
   where $x_{i+1}^{\eps_1} x_i^{\eps_2} x_{i+1}^{\eps_3}\tto V$
   is one of $(iv)$--$(vi)$ and $W$ is a word in
   $ x_1^{\pm1},\dots, x_{i-1}^{\pm1}$;
   \smallskip
\item "$(viii)$"
   $ x_j^{\eps_1} x_i^{\eps_2}\tto x_i^{\eps_2} x_j^{\eps_1}$,
   $\,\{\eps_1,\eps_2\}\subset\{-1,1\}$,  $\,j-1>i\ge1$;
\endroster

   \smallskip
An {\it elementary $K$-reduction of a monomial} is $AUB\to AVB$ where
$AUB\in F_\infty^+$
and $U\to V$ is a basic replacement.
An {\it elementary $K$-reduction of an element of} $kF_\infty^+$ is
$\sum_{j=1}^m c_jW_j \to c_1 W'_1 + \sum_{j=2}^m c_j W_j$ where $c_1,\dots,c_m\in k$,
$W_1,\dots W_m$ are {\sl pairwise distinct}
elements of $F_\infty^+$, and $W_1\to W'_1$
is an elementary $K$-reduction of a monomial.

An element of $F_\infty^+$ (resp. of $kF_\infty^+$)
is {\it $K$-reduced} if no
$K$-reduction can be applied to it. We denote the set of such elements
by $F_\infty^{\red}$ (resp. $kF_\infty^{\red}$). We set also
$F_n^{\red}=F_n^+\cap F_\infty^{\red}$ and
$kF_n^{\red}=kF_n^+\cap kF_\infty^{\red}$. Then $kF_\infty^{\red}$ is a
submodule (not a subalgebra) of $kF^+_\infty$.
We denote $\pi:kF_\infty^+\to K_\infty$ and $\pi_n:kF_n^+\to K_n$
the morphisms of $k$-algebras induced by $x_i\mapsto s_i$.

We say that an element $X$ of $F_\infty^+$ is {\it almost $K$-reduced}
if there exists a sequence $X=X_1\to X_2\to\dots\to X_m$
of elementary $K$-reductions of type $(viii)$ such that $X_m$ is $K$-reduced.

For $X=x_{i_1}^{\eps_1}\dots x_{i_m}^{\eps_m}\in F_\infty^+$, $\eps_j=\pm1$,
we define the {\it weight} $\wt X=\sum_j i_j$
and the {\it auxiliary weight} $\wt' X=\sum_j ji_j$.
It is clear that the set of all monomials of a given weight is finite.
For $X\in kF_\infty^+$ we set
$\wt X=\max_i\wt X_i$ if $X=\sum_i c_iX_i$ with $c_i\in k$ and
$X_1,X_2,\dots$ pairwise distinct elements of $F_\infty^+$.

The following statement is easy and we omit its proof.

\nopagebreak
\proclaim{ Proposition \propReduc }

a). If $X\to X'$ is an elementary $K$-reduction, then $\pi(X)=\pi(X')$ and
%
$\wt X\ge\wt X'$. If, moreover, $X$ is a monomial, then
$\wt X=\wt X'$ if and only if $X\to X'$ is a $K$-reduction of type $(viii)$
and in this case we have $\wt'(X)<\wt'(X')$.

\smallskip
b). $\pi(F_\infty^{\red})$ generates $K_\infty$ as a $k$-module.

\smallskip
c). $kF_\infty^{\red}$ is a free $k$-module and $F_\infty^{\red}$ is a
free base of $kF_\infty^{\red}$.

\smallskip
d). $F_\infty^{\red}$ is the set of all words 
$X_1X_2\dots X_m$ where
   $X_\nu = x_{i_\nu}^{\pm1}x_{i_\nu-1}^{\pm1}\dots x_{j_\nu}^{\pm1}$,
    \hbox{$i_\nu\ge j_\nu$}
    $(1\le\nu\le m)$,
    $i_1<\dots < i_m$, and all the signs are mutually independent.
%

\smallskip
e). {\rm(}Proven in {\rm[\refFunar])} $\pi_3$ is an isomorphism of $k$-modules
$kF_3^{\red}$ and $K_3$.
\endproclaim

\noindent{\bf Remark \remOne.}
Let
$$
  S_{i,j}=\{x_i^{\pm1} x_{i-1}^{\pm1}\dots x_j^{\pm1}\}\qquad
  \text{and}\quad
  S_i=\{1\}\cup S_{i,i}\cup S_{i,i-1}\cup\dots\cup S_{i,1}.
                                                                \eqno(\eqDefS)
$$
Then Part (d) of Proposition \propReduc\
can be stated as follows: each element of $F_n^{\red}$ can be
represented in a unique way as a product $X_1X_2\dots X_{n-1}$ with $X_i\in S_i$.
Since $|S_i|=1+2+\dots+2^i=2^{i+1}-1$, we obtain $|F_n^{\red}|=\prod_{i=1}^n(2^i-1)$,
in particular,
$$
  |F_2^{\red}|=3,\;\;
  |F_3^{\red}|=3\cdot7=21,\;\;
  |F_4^{\red}|=3\cdot7\cdot 15=315,\;\;
  |F_5^{\red}|=3\cdot7\cdot 15\cdot 31=9765.
$$

\noindent{\bf Remark \remTwo.}
In basic replacements $(vii)$, it is enough to consider only words $W$ belonging to
$S_{i-1}$ (see (\eqDefS) for the definition of $S_{i-1}$).

\medskip

We define a $k$-linear mapping $\rr:kF_\infty^+\to kF_\infty^{\red}$ as follows.
For each $X\in F_\infty^+$ we fix an arbitrary sequence of elementary $K$-reductions
$X=X_1\to X_2\to\dots\to X_m\in kF_\infty^{\red}$ and we set $\rr(X)=X_m$.
Then we extend the mapping to $kF_\infty^+$ by linearity.


\subhead\sectMarkov. Markov trace
\endsubhead
Let $A=k[u,v]$ and $AK_n=K_n(\alpha,\beta;A)$. Let $M=M(\alpha,\beta;k)$
be the quotient of $AK_\infty$ by the relations
(\eqMT) and let $t:AK_\infty\to M$ be the quotient map.
We call $t$ {\it the universal Markov trace} on $K_\infty$ over $k$.
It is indeed universal in the sense that any Markov trace on $K_\infty(\alpha,\beta;A)$
with values in an $A$-module $M'$ is $f\circ t$ for some
$f\in\Hom_A(M,M')$.

We define $A$-linear mappings
$\tau_n:AF_n^+\to AF_{n-1}^{\red}$ called {\it Markov reductions}
as follows. By Proposition \propReduc(d), we have
$F_n^{\red}\subset F_{n-1}^{\red}\cup (F_{n-1}^{\red}x_{n-1} F_{n-1}^{\red})\cup
(F_{n-1}^{\red}x_{n-1}^{-1}F_{n-1}^{\red})$.
So, we set $\tau_n(X)=X$, $\tau_n(Xx_{n-1}Y)=u\rr(XY)$, and
$\tau_n(Xx_{n-1}^{-1}Y)=v\rr(XY)$ for $X,Y\in F_{n-1}^+$ and then we extend $\tau_n$
to $AF_n^{\red}$ by linearity and to $AF_n^+$ by setting
$\tau_n(X)=\tau_n(\rr(X))$.
Finally, we define $\tau:F_\infty^+\to AF_1^+=A$
by setting $\tau(X)=\tau_2\circ\dots\circ\tau_n(X)$ for $X\in AF_n^+$.

By definition of $t$ and $\tau$, we have $t(\pi(X))=t(\pi(\tau(X)))$, thus $M=t(K_\infty)$
is generated by $t(1)$.
Let $I=I(\alpha,\beta;k)$ be the annihilator of $M$.
Thus we have $M\cong A/I$.


\subhead \sectMainResult. Statement of the main result
\endsubhead
Let $\sh^n:AF_\infty^+\to AF_\infty^+$, $n\in\Z$, be the $A$-algebra
endomorphism (the $n$-{\it shift\/}) induced by
$$
    \sh^n x_i = \cases x_{i+n}, &i+n>0,\\
                       0,       &i+n\le 0.
		\endcases
$$
We set $\sh=\sh^1$.

For $X\in F_5^+$, we define $\rho_X\in\End_A(AF_4^{\red})$ by setting
$\rho_X(Y)=\tau_5(X \sh Y)$.
Let $J_4=J_4(\alpha,\beta;k)$ be the minimal submodule of $AF_4^{\red}$
satisfying the following properties (recall that the sets $S_{i,j}$ and $S_i$
are defined in (\eqDefS)):
\roster
\item "(J1)"
       $\rr(\rr(X_3 X_2)X_1)-\rr(X_3\rr(X_2X_1))\in J_4$
       for any $X_j\in\sh^{3-j}S_j\setminus\{1\},\,j=1,2,3$;
\smallskip
\item "(J2)"
       $\rho_X(J_4)\subset J_4$ for any $X\in S_4$.
\endroster

\smallskip
In a similar way we define a module $L$. Let
$N=AF_2^{\red}\otimes_A AF_2^{\red}$. We define $A$-linear mappings
$\tau_N:N\to A$ and $\rho_\delta:N\to N$, $\delta=(\delta_1,\delta_2)\in\{-1,0,1\}^2$,
by setting for any $Y=x_1^{\eps_1}\otimes x_1^{\eps_2}$ ($\eps_1,\eps_2\in\{-1,0,1\}$)
$$
  \tau_N(Y)=\tau(x_1^{\eps_1}x_1^{\eps_2}),\qquad
  \rho_\delta(Y)=
      x_1^{\delta_1}\otimes\tau_3(x_2^{\eps_1}x_1^{\delta_2}x_2^{\eps_2})
$$
and we define $L$ as the minimal submodule of
$N$ satisfying the conditions:

\roster
\item "(L1)"
      $\tau_3(x_2^{\eps_1}x_1^{\eps_2}x_2^{\eps_3})\otimes x_1^{\eps_4} - 
      x_1^{\eps_2} \otimes\tau_3(x_2^{\eps_3}x_1^{\eps_4}x_2^{\eps_1})\in L$
      for any $\eps_1,\eps_3\in\{-1,1\}$ and for any $\eps_2,\eps_4\in\{-1,0,1\}$;
\smallskip
\item "(L2)"
      $\rho_\delta(L)\subset L$
      for any $\delta\in\{-1,0,1\}^2$.
\endroster

\proclaim{ Theorem \thMain. (Main Theorem) } $I=\tau(J_4)+\tau_N(L)$.
\endproclaim

The theorem is proved in \S3 (see Proposition \propEasy\ for the inclusion ``$\supset$''
and Proposition \propDiffi\ for the reverse inclusion). 


This result allows (at least theoretically) to compute $I$. Indeed,
we start with the $A$-module $J_4^{(0)}$ generated by the elements in (J1)
and and we set $J_4^{(i+1)}=\sum_{X\in S_4}\rho_X(J_4^{(i)})$.
Then the Gr\"obner bases $G^{(i)}$ of the modules $J_4^{(i)}$ can be
computed recursively using the fact that $J_4^{(i+1)}$ is the module
generated by $\bigcup_{X\in S_4}\rho_X(G^{(i)})$. So, we construct
an increasing sequence of submodules $J_4^{(0)}\subset J_4^{(1)}\subset\dots$.
Since the subring of $A$ generated by $\alpha,\beta,u,v$
is noetherian, there exists $m_0$ such that $J_4^{(m_0)}=J_4^{(m_0+1)}=\dots$
($m_0$ is determined by the condition $G^{(m_0)}=G^{(m_0+1)}$).
The module $L$ can be computed in a similar way as the limit of
$L^{(0)}\subset L^{(1)}\subset\dots$ where $L^{(0)}$ is generated by the elements in (L1)
and $L^{(i+1)}=\sum_\delta \rho_\delta(L^{(i)})$.

Performing in practice this computation
for $\alpha=\beta=0$, $k=\Z$
(the case considered in [\refFunar] and [\refCM]) and in some other special cases,
we obtain the following results.
To compute the Gr\"obner bases, we used {\tt Singular 3-1-3} and {\tt Macaulay2} software;
see details on the web page [\refWebPage].

\proclaim{ Corollary \corZZZ }
$I(0,0;\Z)=(16,\,4u^2+4v,\,4v^2+4u,\,u^3+v^3+uv-3)$.
\endproclaim

\proclaim{ Corollary \corBetaZero } Let $k=\bold k[\alpha]$ for a ring $\bold k$ specified below and
let $I=I(\alpha,0,k)$. Let $\Cal G$ be the reduced Gr\"obner base of $I$
with respect to the lexicographic order such that $v>u>\alpha$ (the Gr\"obner bases in (e) -- (g)
also correspond to this order).
Let
$$
\split
   &f_1=\gamma_1\gamma_2\gamma_3,\qquad\text{where }\;
   \gamma_1=\alpha^3+8,\;
   \gamma_2=2\alpha^3+1,\;
   \gamma_3=3\alpha^3+8,
\\&
   f_2=\gamma_1\gamma_3(u-\alpha),
\\&
   f_3=\gamma_3(6u^3 - 3\alpha^2 u + \alpha^3+2),
\\&
   f_4=336u^4-792\alpha u^3
   +12(15\alpha^3+106)\alpha^2u^2 
   +6(141\alpha^3+544)u 
\\&\qquad
   -114\alpha^7-1405\alpha^4-3152\alpha,
\\&
   f_5=288v+336\alpha^2u^3 
   +72(3\alpha^3+28)u^2  
   -48(9\alpha^3+44)\alpha u 
   -6a^8+53a^5+472a^2.
\endsplit
$$

\smallskip
a). If $\bold k=\Bbb F_2$, then
$\Cal G=\{\alpha^4,\,
\alpha^2(u^3+1+\alpha u^2),\,
\alpha^2(v+u^2),\,
v^3+u^3+uv+1+\alpha(uv^2+v+\alpha u)\}$, hence $\dim_{\bold k}(A/I)=\infty$.
\smallskip



b). If $\bold k=\Bbb F_3$, then
$\Cal G=\{
\alpha^3-1,\,
(u^2-\alpha^2)(u^2-\alpha u-\alpha^2),\,
v+u^2
\}$, hence $\dim_{\bold k}(A/I)=12$.
\smallskip

c). If $\bold k=\Bbb Q$ or $\bold k=\Bbb F_p$ for a prime $p$ in the range $5\le p\le 397$, $p\ne37$
{\rm(conjecturally, 
for any prime $p\not\in\{2,3,37\}$)},
then $I=(f_1,\dots,f_5)$. We have $\Cal G=\{f_1,\dots,f_5\}$ and hence $\dim_{\bold k}(A/I)=24$
except the case $\bold k=\Bbb F_7$ where we have 
%
$\Cal G=\{f_1,f_2,
u^3+2\alpha u(\gamma_1u+\gamma_3\alpha)+3\alpha^6-3\alpha^3,f_5\}$
and hence $\dim_{\bold k}(A/I)=21$.

\smallskip
d). If $\bold k=\Bbb F_{37}$, then we have $I=(f_1,f_2,f_3,(u+7\alpha)f_4,f_5-14\alpha f_4)$ and
$\Cal G=\{f_1,f_2,f_3, f_{4,37},f_{5,37}\}$ where
$$
\split
&f_{4,37}=u^5 + 2\alpha u^4 + 7\alpha^2 u^3 - 9(\alpha^3-1)u^2 + (6\alpha^3+2)\alpha u-12a^8+9a^5, \\
&f_{5,37}=v-4\alpha u^4 + 15\alpha^2u^3 - (14\alpha^3+16)u^2 - (\alpha^3+18)\alpha u -4\alpha^8+2\alpha^5-2\alpha^2
\endsplit
$$
and hence $\dim_{\bold k}(A/I)=27$.

\smallskip
e).
If $\bold k=\Bbb Z/32\Bbb Z$, then
$\{8\alpha^3,\,4\alpha^4,\,\alpha^5+8\alpha^2,\,16\alpha u^3+16\alpha^2u^2+4\alpha^3u+2\alpha^4,\, 
   8\alpha^2u^3-12\alpha^3u^2+8\alpha^2,\,2\alpha^3u^3+16u^3+16\alpha^2u+6\alpha^3+16,\,           
   \alpha^4u^3+8\alpha u^3+4\alpha^5u^2+16\alpha^2u^2+\alpha^4+8\alpha,\,                          
   16u^4-8\alpha u^3+8\alpha^2u^2-6\alpha^3u+\alpha^4+16\alpha,\,                                  
   4\alpha^2u^4+7\alpha^3u^3-8u^3-2\alpha^4u^2-4\alpha^2u-9\alpha^3-8,\,                           
   8\alpha u^5+16u^3+\alpha^4u^2+8\alpha^2u+6\alpha^3,\,                                           
   \alpha^3u^5-8u^5-8\alpha u^4-4\alpha^2u^3-5\alpha^3u^2+8u^2-13\alpha^4u+16\alpha u-4\alpha^2,\, 
   4\alpha u^6-9\alpha^3u^4+8u^4+3\alpha^3u+8u-3\alpha^4+4\alpha,\,                                
   16v+8\alpha u^4+12\alpha^2u^3+16u^2+3\alpha^4u-12\alpha^2,\,                                    
   8\alpha v+10\alpha^4u^2-8\alpha u^2+16\alpha^2u+2\alpha^3+16,\,                                 
   4\alpha^2v-12\alpha^2u^2+12\alpha^3u-13\alpha^4-8\alpha,\,  \alpha^3v-8v+\alpha^3u^2-8u^2,\,    
   2\alpha^2uv+8v+4\alpha u^4-12\alpha^2u^3-4\alpha^3u^2-8u^2+13\alpha^4u+12\alpha u-2\alpha^2,\,  
  4\alpha u^3v+12\alpha v-12\alpha u^5+13\alpha^3u^3+8u^3-10\alpha^4u^2+12\alpha u^2+3\alpha^3-8,\,
  4v^2-2\alpha uv+3\alpha^2v+\alpha^2u^2+4u-2\alpha,\,                                             
  2\alpha v^2+15\alpha^2uv-4v+15\alpha^3u^2-4u^2+2\alpha u+15\alpha^2,\,                           
  \alpha^2v^2+4uv-2\alpha v-4\alpha u^5+15\alpha^4u^2-2\alpha u^2-7\alpha^2u+13a^3-4,\,            
  v^3-15\alpha uv^2-3uv-7\alpha v+u^3+6\alpha u^2+\alpha^2u+2\alpha^3-15\}$                        
is a Gr\"obner base of $I$.
\if01{
If $\bold k=\Bbb Z/16\Bbb Z$, then
$\{4\alpha^3,\,
2\alpha^4,\,
\alpha^5+8\alpha^2,\,
8\alpha u^3+8\alpha^2 u^2+2\alpha^3 u+\alpha^4,\,
4\alpha^2 u^3+2\alpha^3 u^2+4\alpha^2,\,
\alpha^3 u^3+8u^3+8\alpha^2 u+3\alpha^3+8,\,
4\alpha u^6+8\alpha^2 u^2+2\alpha^3 u+\alpha^4+4\alpha,\,
8\alpha v+8\alpha u^2+2\alpha^3,\,
4\alpha^2 v+4\alpha^2 u^2+\alpha^4+8\alpha,\,  
\alpha^3 v+8v+\alpha^3 u^2+8u^2,\,    
2\alpha^2 uv+8v-4\alpha u^4+8u^2-4\alpha u+2\alpha^2,\,  
4\alpha u^3v+4\alpha v-4\alpha u^5+\alpha^4 u^2+4\alpha u^2,\,
2\alpha v^2+7\alpha^2 uv-7\alpha^3v+4v+4u^2+2\alpha u+7\alpha^2,\,
\alpha^2v^2+4uv-2\alpha v+4\alpha u^5-2\alpha u^2+\alpha^2 u+3\alpha^3-4,\,
v^3+\alpha uv^2-3uv-7\alpha v+u^3+6\alpha u^2+\alpha^2u+2\alpha^3+1\}$
is a Gr\"obner base of $I$.
$\Cal G=\{4\alpha^3,\,
2\alpha^4,\,
\alpha^2\gamma_1,\,
8\alpha u^3+8\alpha^2 u^2+2\alpha^3 u+\alpha^4,\,
\alpha^2(4u^3+2\alpha u^2+4),\,
\gamma_1(u^3+3) + 8\alpha^2 u,\,
4\alpha u^6+8\alpha^2 u^2+2\alpha^3 u+\alpha^4+4\alpha,\,
8\alpha(v+u^2)+2\alpha^3,\,
4\alpha^2(v+u^2)+\alpha\gamma_1,\,  
\gamma_1(v+u^2),\,    
2\alpha^2(uv+1)+8(v+u^2)+12\alpha u(u^3+1),\,  
4\alpha v(u^3+1)+12\alpha u^5+\alpha^4 u^2+4\alpha u^2,\,
2\alpha(v^2+u)+7\alpha^2(uv+1)+9\alpha^3 v+4(v+u^2),\,
\alpha^2(v^2+u)+4(uv-1)-2\alpha(v+u^2)+4\alpha u^5+3\alpha^3,\,
v^3+u^3-3uv+1+\alpha(uv^2+9v+6u^2+\alpha u+2\alpha^2)\}$.
}\fi 

\smallskip
f). If $\bold k=\Bbb Z/3^r\Bbb Z$ for $r\le6$ {\rm(conjecturally, 
for any $r$)}, then
$I=(\Cal F)$ for $\Cal F=\{f_1,f_2,f_3,\,2u^4 + \alpha u^3-\alpha^2 u^2 + 2u - \alpha,\, v+u^2\}$.
If, moreover, $r\ge 2$, then $\Cal F\cup\{3^{r-1}(\alpha^3-1)\}$
is a Gr\"obner base of $I$,
hence $|A/I|=3^{12r}$. By (b), this implies
that $A/I$ is a free $\bold k$-module of rank $12$.

\smallskip
g). If $\bold k=\Bbb Z/37p\Bbb Z$ for a prime $p\le 67$, $p\not\in\{2,3,7\}$
{\rm(conjecturally, for any prime $p\not\in\{2,3,7\}$)}, then $\{f_1,f_2,f_3,37 f_4, (u+7\alpha)f_4,
f_5-14\alpha f_4\}$ is a Gr\"obner base of $I$.
\endproclaim
\if01{
//p=2
ring r=2,(v,u,a),lp;
poly h1=a4;
poly h2=a2u2+a2v;
poly h3=a2v2+a2u+a3;
poly h4=u3+auv2+v3+uv+a2u+av+1;
ideal I=h1,h2,h3,h4;
I=std(I);I;

I[1]=a4
I[2]=u3a2+u2a3+a2
I[3]=va2+u2a2
I[4]=v3+v2ua+vu+va+u3+ua2+1
>

ring r=2,(V,U,v,u,a),(lp(2),lp);
poly h1=a4;
poly h2=a2u2+a2v;
poly h3=a2v2+a2u+a3;
poly h4=u3+auv2+v3+uv+a2u+av+1;
ideal J=h1,h2,h3,h4,Uu-1,Vv-1;
J=std(J);J;

J[1]=a4
J[2]=u3a2+u2a3+a2
J[3]=va2+u2a2
J[4]=v3+v2ua+vu+va+u3+ua2+1
J[5]=Ua2+v2ua2+vu2a3+ua3
J[6]=Uu+1
J[7]=Va2+U2a2
J[8]=Vu3+Vua2+V+v2+vua+u+a
J[9]=Vv+1
J[10]=VU+Vu2+Va2+Uv2+Uvua+Uu+Ua

//p=3
ring r=3,(v,u,a),(lp(2),lp);
poly h1=a3-1;
poly h2=(u2-a2)*(u2-au-a2);
poly h3=v+u2;
ideal I=h1,h2,h3; I=std(I);I;

ring r=3,(V,U,v,u,a),(lp(2),lp);
poly h1=a3-1;
poly h2=(u2-a2)*(u2-au-a2);
poly h3=v+u2;
ideal J=h1,h2,h3,Uu-1,Vv-1; J=std(J);J;

ring r=11,(v,u,a),lp;
poly g(1)= a3+8;
poly g(2)=2a3+1;
poly g(3)=3a3+8;
poly f(1)=g(1)*g(2)*g(3);
poly f(2)=g(1)*g(3)*(u-a);
poly f(3)=g(3)*(6u3-3a2u+a3+2);
poly f(4)=336u4-792u3a+180u2a5+1272u2a2+846ua3+3264u-114a7-1405a4-3152a;
poly f(5)=288v+336u3a2+216u2a3+2016u2-432ua4-2112ua-6a8+53a5+472a2;
ideal Id; Id=f(1),f(2),f(3),f(4),f(5); Id=std(Id); Id;
reduce(Id,II);
reduce(II,Id);

ring r=11,(V,U,v,u,a),(lp(2),lp);
poly g(1)= a3+8;
poly g(2)=2a3+1;
poly g(3)=3a3+8;
poly f(1)=g(1)*g(2)*g(3);
poly f(2)=g(1)*g(3)*(u-a);
poly f(3)=g(3)*(6u3-3a2u+a3+2);
poly f(4)=336u4-792u3a+180u2a5+1272u2a2+846ua3+3264u-114a7-1405a4-3152a;
poly f(5)=288v+336u3a2+216u2a3+2016u2-432ua4-2112ua-6a8+53a5+472a2;
ideal Id1; Id1=f(1),f(2),f(3),f(4),f(5),Uu-1,Vv-1; Id1=std(Id1); Id1;

ideal id1;ideal id2;ideal id3;
id1=Id,g(1); id1=std(id1); 
id2=Id,g(2); id2=std(id2); 
ideal tmp; tmp=id2[1];
id2[2]=reduce(id2[2],tmp);
tmp=tmp,id2[2];
id2[3]=reduce(id2[3],tmp);
id3=Id,g(3); id3=std(id3); 
ideal tmp; tmp=id3[1];
id3[2]=reduce(id3[2],tmp);
tmp=tmp,id3[2];
id3[3]=reduce(id3[3],tmp);
14u^4-33u^3a+33u^2a^2+42u*(-3/8*a^3)-9a*(-3/8*a^3)

336u^4-792u^3a+180u^2a^5+1272u^2a^2+846u a^3+3264u-114a^7-1405a^4-3152a

//p=7

I[1]=a9+3a6+a3-1
I[2]=ua6-ua3-2u-a7+a4+2a
I[3]=u3+2u2a4+2u2a-ua5+2ua2+3a6-3a3
I[4]=v+3u4a-u2+2ua4-ua-3a5-3a2

ideal I;
I=
a9+3a6+a3-1,
ua6-ua3-2u-a7+a4+2a,
u3+2u2a4+2u2a-ua5+2ua2+3a6-3a3,
v+3u4a-u2+2ua4-ua-3a5-3a2;
}\fi

\noindent{\bf Remark \remInv. }
The Markov trace $t$ over $k$ defines an invariant of oriented links
$P(L)=P_{\alpha,\beta,k}(L) = u^{(1-n-e)/2} v^{(1-n+e)/2}t(b)\in
k[u^{\pm1/2},v^{\pm1/2}]/I(\alpha,\beta;k)$
where $b$ is a representation of
a link $L$ by a braid with $n$ strings and $e$ is the sum of exponents of $b$.
It is shown in [\refCM] that $P_{0,0;\Z/4\Z}$
is determined by the HOMFLY polynomial (see Section \sectHOMFLY\ below).
It is not known whether this is true or not in other cases.
In the first arxiv version of this paper (arxiv:1206.0765v1) it was claimed
that $P_{\alpha,0;\Bbb Q[\alpha]}$ and $P_{\alpha,0;\Bbb F_3[\alpha]}$
detect the chirality of the knot $10_{71}$. Unfortunately, this is not so.
However, $P_{\alpha,0;\Bbb F_{37}[\alpha]}$
detect the chirality of the knots $10_{48}$, $10_{91}$ and it distinguishes many other pairs
of knots with equal HOMFLY polynomials.
In the cases computed so far, the invariants $P_{\alpha,0,k}$ do not distinguish
any pair of knots up to 11 crossings with equal Kauffman polynomials.


\medskip

\noindent{\bf Remark \remMarin. }
According to [\refMarin; Theorem 1.5], we have $I(\alpha,\beta;\Bbb Q(\alpha,\beta))=(1)$.
This is equivalent to say that $I(\alpha,\beta;\Bbb Q[\alpha,\beta])$ contains
a nonzero polynomial in $\alpha$ and $\beta$
(observe that the same phenomenon takes place in Corollaries \corZZZ\ and \corBetaZero).

\medskip

%

%
\noindent{\bf Remark \remTransv. } 
Our main motivation for studying
quotients of cubic Hecke algebras and their Markov traces stems from the
Markov-type theorem for transversal links proved in [\refOS]. If $u, v$ were zero divisors in
$A/I$ then one might possibly obtain transversal links invariants which distinguish
isotopic links with the same Bennequin invariants. In the cases computed so far,
however, $u$ and $v$ are not zero divisors.

We say that $t\in\Hom_{k[u]}(K_\infty\otimes k[u],M)$ is a {\it semi-Markov trace} on
$K_\infty$ if $t(xy)=t(yx)$ and $t(xs_n)=ut(x)$ for $x,y\in K_n$, $n>0$.
Due to [\refOS], any semi-Markov trace provides an invariant of transversal links.
In [\refO], the methods of the present paper are adapted for studying the universal semi-Markov
trace on $K_n$.
\medskip

\noindent{\bf Remark \remFunarErr.} The main error 
in [\refFunar] (which was
repeated also in [\refBF]) is that the modules $J_4^{(0)}$ and $L^{(0)}$
were considered instead of $J_4$ and $L$.


\head 3. Proof of Main Theorem
\endhead

The idea of the proof is as follows.
Given $X\in AF^+_n$, the element of $A$ representing $\tau(X)$ is computed by successive reductions (i)--(viii),
Markov relations, and cyclic permutations. We have to find the minimal possible ideal $I$ such
that the result does not depend on the order of these operations.
It is easy to observe (though not so easy to formalize this observation)
that the main sources of the ambiguity are as follows. First,
the reduction of a subword
$W=x_3^{\eps_1} x_2^{\eps_2} x_1^{\eps_3} x_3^{\eps_4} x_2^{\eps_5} x_3^{\eps_6}$
can be started either with $x_3^{\eps_4} x_2^{\eps_5} x_3^{\eps_6}$ or
(after commuting $x_1^{\eps_3}$ and $x_3^{\eps_4}$) with
$x_3^{\eps_1} x_2^{\eps_2} x_3^{\eps_4}$.
Thus $\tau(J_4^{(0)})$ should be included in $I$.
Second, the reduction of a word $x_2^{\eps_1} x_1^{\eps_2} x_2^{\eps_3} x_1^{\eps_4}$ can be started
either by $x_2^{\eps_1} x_1^{\eps_2} x_2^{\eps_3}$ or (after a cyclic permutation) by
$x_2^{\eps_3} x_1^{\eps_4} x_2^{\eps_1}$. Thus $\tau_N(L^{(0)})$ should be included into $I$.

Let us focus on the first case.
So, let $Y$ be an element of $J_4^{(0)}$. Then, for any $X$ and any $n$,
we should have $\tau(X\sh^n Y)\in I$. In particular, for any $X\in F_5^+$, we should have
$\tau(\rho_X(Y))\in I$, i.~e., $\tau(J_4^{(1)})\subset I$. By iterating this process, we
conclude that $\tau(J_4)\subset I$. Similarly, $\tau_N(L)\subset I$. This is the easy part of the proof
which is formally exposed in Section 3.1.

The difficult part of the proof (formally exposed in Section 3.2) consists in checking that
any choice of the reduction process leads to the same result modulo $I'=\tau(J_4)+\tau_N(L)$.
We use induction on the weight (see the definition of the weight function wt in Section 2.1).
As we pointed out above, there are two main sources of the ambiguity. Again, we discuss here only the
first one. So, we have to prove that $\tau(X\sh^n Y)\in I'$ for any $X$ when $Y\in J_4$.
By additivity, we may assume that $X$ is a monomial.
If any reduction can be applied to $X$, then we do it and we use the induction hypothesis.
So, we may assume that $X\in AF_{n+4}^+$. 
If $X=X_1X_2$ where $X_2$ commutes with $\sh^{n} Y$, then we replace $X\sh^n Y$ by
$X_2X_1\sh^n Y$. Thus we arrive to the case when $X=X'\sh^{n-1} X_1$ with $X'\in F_{n+3}^+$, $X_1\in F_5^+$,
and we apply the induction hypothesis to $X'\sh^{n-1}Y'$ where $Y' = \rho_{X_1}(Y)$.

\subhead 3.1. Easy part: $\tau(J_4)+\tau_N(L)\subset I$
\endsubhead

Let $J_4^{(0)}\subset J_4^{(1)}\subset\dots$ and $L^{(0)}\subset L^{(1)}\subset\dots$
be as defined in \S\sectMainResult.


For $n\ge 4$ and $a\in AK_n$, we define
$t_{n,a}\in\Hom_A(F_4^{\red}, A)$ by setting
$t_{n,a}(X)=t(a\,\pi(\sh^{n-4}X))$. Similarly,
for $n\ge 1$ and $a,b\in AK_n$, we define
$t_{n,a,b}\in\Hom_A(N,A)$ by setting
$t_{n,a,b}(X\otimes Y)=t\big(\pi(\sh^{n-1}X)\,a\,\pi(\sh^{n-1}Y)\,b\big)$.

\proclaim{ Lemma \lemEasy } 

a). $J_4\subset\ker t_{n,a}$ for any $n\ge4$ and any $a\in K_n$.

b). $L\subset\ker t_{n,a,b}$ for any $n\ge1$ and any $a,b\in K_n$.
\endproclaim

\demo{ Proof }
We prove by induction that
a) $J_4^{(i)}\subset\ker t_{n,a}$ and b)
$L^{(i)}\subset\ker t_{n,a,b}$. For $i=0$, the statement is evident.
Suppose that it is true for $i-1$ and let us prove it for $i$.
Note that we have
$$
    t\big(a\,\pi(\sh^p\tau_{n-p}(X))\,b\big)=t\big(a\,\pi(\sh^p X)\,b\big)
   \qquad\text{for}\quad
    a,b\in K_{n-1}, \; X\in AF_n^+
                                                   \eqno(\eqTau)
$$
\smallskip
a). It is enough to check that $\rho_X(Y)\in\ker t_{n,a}$
for any $Y\in J_4^{(i-1)}$,
$X\in S_4$. $n\ge 4$, $a\in K_n$. Indeed,
$$
\xalignat2
  t_{n,a}(\rho_X(Y))&=
  t\big(a\,\pi(\sh^{n-4}\rho_X(Y))\big)
                     &&\text{by definition of $t_{n,a}$}
\\&
  =t\big(a\,\pi(\sh^{n-4} \tau_5(X \sh Y) )\big)
                     &&\text{by definition of $\rho_X$}
\\&
  =t\Big(a\,\pi\big((\sh^{n-4}X)(\sh^{n-3}Y)\big)\Big)
                     &&\text{by (\eqTau)}
\\&
  =t_{n+1,a'}(Y)
                     &&\text{for $a'=a\,\pi(\sh^{n-4}X)\in K_{n+1}$}
\\&
  =0                 &&\text{by the induction hypothesis}
\endxalignat
$$

\smallskip
b). It is enough to check that $\rho_\delta(Y)\in\ker t_{n,a,b}$
for any $Y\in L^{(i-1)}$, $\delta=(\delta_1,\delta_2)\in\{-1,0,1\}^2$,
$n\ge 1$, $a,b\in K_n$. Indeed, let
$Y=\sum_j c_j x_1^{\eps_1(j)}\otimes x_1^{\eps_2(j)}$. Then
$$
\xalignat2
  t_{n,a,b}&(\rho_\delta(Y))
  =t_{n,a,b}\Big(\sum c_j
  x_1^{\delta_1}\otimes\tau_3(x_2^{\eps_1(j)}x_1^{\delta_2}x_2^{\eps_2(j)})
  \Big)
        &&\text{def. of $\rho_\delta$}
\\&
  =\sum c_j t\Big(\pi(\sh^{n-1}x_1^{\delta_1})\,a\,
    \pi\big(\sh^{n-1}\tau_3(x_2^{\eps_1(j)}x_1^{\delta_2}x_2^{\eps_2(j)}
    )\big)\,b\Big)\;
        &&\text{def. of $t_{n,a,b}$}
\\&
  =\sum c_j t\big(s_n^{\delta_1}\,a\,
    s_{n+1}^{\eps_1(j)}s_n^{\delta_2}s_{n+1}^{\eps_2(j)}
    \,b\big)
                     &&\text{by (\eqTau)}
\\&
  =\sum c_j t\big(
    s_{n+1}^{\eps_1(j)}s_n^{\delta_2}s_{n+1}^{\eps_2(j)}
    \,b\,s_n^{\delta_1}\,a\big)
                     &&t(xy)=t(yx)
\\&
  =\sum c_j t\Big(
    \pi\big(\sh^n x_1^{\eps_1(j)}\big)\,s_n^{\delta_2}
    \pi\big(\sh^n x_1^{\eps_2(j)}\big)\,
    \,b\,s_n^{\delta_1}\,a\Big)
\\&
  =t_{n+1,a',b'}(Y)
                     &&\hskip-12mm
		     a'=s_n^{\delta_2},\,b'=bs_n^{\delta_1}a
\\&
  =0                 &&\hskip-12mm
                    \text{by induction hypothesis}	     
\endxalignat
$$
\enddemo

\proclaim{ Proposition \propEasy } $\tau(J_4)+\tau_N(L)\subset I$.
\endproclaim

\demo{Proof}
Indeed, by Lemma \lemEasy, we have
$t(\tau(X))=t_{4,1}(X)=0$ for any $X\in J_4$ and
$t(\tau_N(X))=t_{1,1,1}(X)=0$ for any $X\in L$.
Thus $\tau(J_4)+\tau_N(L)\subset\ker (t|_A)=I$.
\enddemo


\subhead 3.2. Difficult part: $I\subset\tau(J_4)+\tau_N(L)$
\endsubhead

Let, as above, $\bar R$ be the submodule of $K_\infty$ generated
by the elements (\eqMT). Set $R=\pi^{-1}(\bar R)$.
Then we have $I=A\cap\bar R=A\cap R$.
Let $\wt: AF_\infty^+\to\Z_{\ge0}$ be the weight function defined in
\S\sectReduc. It defines a filtration on $AF_\infty^+$, namely,
$A=AF_{[0]}^+\subset AF_{[1]}^+\subset AF_{[2]}^+\subset\dots$
where $AF_{[w]}^+=\{X\in AF_\infty^+\, |\, \wt X\le w\}$.

We shall work with the following set
of generators
$\Cal R=\Cal R_T\cup \Cal R_M\cup \Cal R_N\cup  \Cal R_H$
of $R$ as an $A$-module (we set here $u_+=u$, $u_-=v$):
$$
\xalignat2
   &\Cal R_T=\{XY-YX\,|\,X,Y\in F_\infty^+\},
      &&\text{trace relations;}
\\&
  \Cal R_M=\{x_n^{\pm1}X-u_\pm X\,|\,X\in F_n^+, n\ge1\},
      &&\text{Markov relations;}
\\&
  \Cal R_N=\{UX-VX\,|\,X,U\in F_\infty^+,
  U\overset{(i)\text{--}(vi)}\to\longrightarrow V\},
      &&\text{nonhomogeneous $K$-relations;}
\\&
  \Cal R_H=\{UX-VX\,|\,X,U\in F_\infty^+,
  U\overset{(viii)}\to\longrightarrow V\},
      &&\text{homogeneous $K$-relations.}
\endxalignat
$$
Let $\Cal R_{[w]}=\Cal R\cap AF_{[w]}^+$, let 
$R_{[w]}$ be the $A$-submodule of $R$ generated by $\Cal R_{[w]}$,
%
%
and let $H$ be the submodule generated by $\Cal R_T\cup\Cal R_H$
(the elements of $H\cap R_{[w]}$ are $\wt$-homogeneous for any $w$).
%
%
Note, that by Proposition \propReduc(a) 
we have
$$
  X\equiv \rr(X)\equiv\tau_n(X)\equiv\tau(X)\mod R_{[\wt X]}
  \qquad\text{for $X\in AF_n^+$.}
                                                     \eqno(\eqWt)
$$
In what follows, a notation like $X_1\equiv X_2\equiv X_3\equiv\dots$
means that $X_i\equiv X_{i+1}\mod R_{[\wt X_i]}$ and $\wt X_{i+1}\ge\wt X_i$,
in particular, in this case we always have
$X_1\equiv X_2\equiv X_3\equiv\dots\mod R_{[\wt X_1]}$.

\proclaim{ Lemma \lemDiffiOne }
Let $Z=X\sh^{n-4} Y$ for $X\in AF_\infty^+$,
$Y\in J_4\cap\sh^{4-n}AF_\infty^+$, $n\ge1$.
Then $Z\in R_{[w]}+\tau(J_4)$ where $w=\wt Z$.
\endproclaim

\demo{ Proof }
We denote $\sh^{n-4} Y$ by $Y_n$. If $X\in AF_m^+$ with $m>n$, then
$$
  XY_n\equiv\tau_m(X)Y_n\equiv\tau_{m-1}(\tau_m(X))Y_n\equiv\dots\equiv
  \tau_{n+1}\circ\dots\circ\tau_{m-1}\circ\tau_m(X)Y_n,
$$
hence it is enough to prove the statement of the lemma under the
additional hypothesis $X\in AF_{n}^+$. We prove it by induction.

If $n=1$, then $X\in AF_1^+=A$ and
$Y\in J_4\cap\sh^3 AF_\infty^+ = J_4\cap A\subset\tau(J_4)$, so,
the statement is trivial.

Suppose that $n\ge 2$, the statement is true for $n-1$, and let us prove it
for $n$. By linearity, it is enough to consider the case when $X\in F_{n}^+$
and since $X\equiv\rr(X)$, we may assume that $X\in F_{n}^{\red}$.
Let $X=X_1X_2\dots X_{n-1}$, $X_i\in S_i$ (see Remark \remOne).
We have $X_{n-1}=(\sh^{n-5}X'_4)X''_{n-5}$ with
$X'_4\in S_4\cap\sh^{5-n}AF_\infty^+$ and
$X''_{n-5}\in S_{n-5}$ (we assume here that $S_i=\{1\}$ when $i\le0$).
Note that $Y_n$ may involve only
$x_{n-4+i}^{\pm1}$, $i=1,2,3$, whereas $X''_{n-5}$
may involve only $x_i^{\pm1}$,
$i\le n-5$, hence they commute.
%
%
Therefore, denoting $X_1\dots X_{n-2}$ by $X'''_{n-2}$, we obtain
$$
\split
  Z&=X'''_{n-2}(\sh^{n-5}X'_4)X''_{n-5}(\sh^{n-4} Y)
   \equiv X'''_{n-2}(\sh^{n-5}X'_4)(\sh^{n-4} Y) X''_{n-5}
\\&
   \equiv X''_{n-5}X'''_{n-2}(\sh^{n-5}X'_4)(\sh^{n-4} Y)
   = X''_{n-5}X'''_{n-2}\sh^{n-5}(X'_4\sh Y)
\\&
   \equiv X''_{n-5}X'''_{n-2}\sh^{n-5}\big(\tau_5(X'_4\sh Y)\big)
   = X' \sh^{n-5} Y'
\endsplit
$$
where $X'=X''_{n-5}X'''_{n-2}\in AF_{n-1}^+$
and $Y'=\tau_5(X'_4\sh Y)=\rho_{X'_4}(Y)\in J_4$.

To complete the proof, it remains to check that
$Y'\in\sh^{5-n} AF_\infty^+$.
Indeed, we have $X'_4\in\sh^{5-n}AF_\infty^+$,
$Y\in \sh^{4-n}AF_\infty^+$, hence $\sh Y\in\sh^{5-n}AF_\infty^+$
and we obtain $X'_4\sh Y\in\sh^{5-n}AF_\infty^+$ whence
$Y'=\tau_5(X'_4\sh Y)\in\sh^{5-n}AF_\infty^+$.
\qed\enddemo

The next lemma is similar.
For $n\ge 1$ and $X_1,X_2\in AF_n^+$ we define
$\varphi_{n,X_1,X_2}\in\Hom_A(N,AF^+_{n+1})$ by setting
$\varphi_{n,X_1,X_2}(Y_1\otimes Y_2)=X_1(\sh^{n-1}Y_1)X_2(\sh^{n-1}Y_2)$.

\proclaim{ Lemma \lemDiffiTwo }
Let $Z=\varphi_{n,X_1,X_2}(Y)$ for $n\ge 1$, $X_1,X_2\in AF^+_n$, $Y\in L$.
Then $Z\in R_{[w]}+\tau_N(L)$ where $w=\wt Z$.
\endproclaim

\demo{ Proof }
It is enough to consider the case when $X_1,X_2\in F_n^{\red}$.
Then there exist $X_i',X_i''\in F_{n-1}^{\red}$ and $\delta_i\in\{-1,0,1\}$
such that $X_i=X'_i x_{n-1}^{\delta_i} X''_i$ $(i=1,2)$.
Let
$$
    Y=\sum_j c_j x_1^{\eps_1(j)}\otimes x_1^{\eps_2(j)}.        \eqno(\eqY)
$$
Then we have
$$
\split
   Z&
   =\sum c_j \,   X_1 x_n^{\eps_1(j)}
                  X_2 x_n^{\eps_2(j)}
   =\sum c_j \,   X'_1 x_{n-1}^{\delta_1} X''_1 x_n^{\eps_1(j)}
                  X'_2 x_{n-1}^{\delta_2} X''_2 x_n^{\eps_2(j)}
\\&
   \equiv\sum c_j \,  X''_2 X'_1 x_{n-1}^{\delta_1} X''_1 X'_2 x_n^{\eps_1(j)}
                          x_{n-1}^{\delta_2} x_n^{\eps_2(j)}
\\&
   =\sum c_j \, X''_2 X'_1 x_{n-1}^{\delta_1} X''_1 X'_2 
                         \sh^{n-2}\big(x_2^{\eps_1(j)}
                          x_1^{\delta_2} x_2^{\eps_2(j)}\big)
\\&
   \equiv\sum c_j \,  X''_2 X'_1 x_{n-1}^{\delta_1} X''_1 X'_2 
                         \sh^{n-2}\tau_3\big(x_2^{\eps_1(j)}
                          x_1^{\delta_2} x_2^{\eps_2(j)}\big)
   = \varphi_{n-1,\bar X_1,\bar X_2}(\bar Y)
\endsplit
$$
where $\bar X_1=X''_2 X'_1$, $\bar X_2=X_1'' X'_2$, $\bar Y=\rho_\delta(Y)$.
So, we have $Z\equiv\bar Z=\varphi_{n-1,\bar X_1,\bar X_2}(\bar Y)$
where $\bar X_1,\bar X_2\in AF_{n-1}^+$, $\bar Y\in L$.

Thus, by induction we reduce the problem to the case $n=1$. In this
case we have $X_1,X_2\in AF_1^+=A$, hence, for $Y$ as in (\eqY), we have
$Z=\varphi_{1,X_1,X_2}(Y)=\sum c_j x_1^{\eps_1(j)}x_1^{\eps_2(j)}$,
hence $Z\equiv\tau_2(Z)=\tau_N(Y)\in \tau_N(L)$.
\qed\enddemo


The next statement can be considered as an
improvement of the Pentagon Lemma from [\refFunar].

\nopagebreak
\proclaim{ Lemma \lemPent\ (Pentagon Lemma) }
Let $Z_1,Z_2\in\Cal R_N\cup\Cal R_M$ and
$Z_1-Z_2\in H+AF_{[w-1]}^+$ where $w=\wt Z_1=\wt Z_2$.
Then $Z_1-Z_2\in H+R_{[w-1]}+\tau(J_4)+\tau_N(L)$.
\endproclaim


\demo{ Proof }
Let $X_i\in F_\infty^+$ be the leading monomial of $Z_i$, $i=1,2$, i.~e.,
$\wt X_i=\wt Z_i$ and $\wt(Z_i-X_i)\le w-1$.
Then $X_1-X_2\in H$, hence there exists a sequence
of words $X_1=W_1,\dots,W_m=X_2$ such that $W_{i+1}$ is obtained
from $W_i$ either by a cyclic permutation or by exchanging two
consecutive commuting letters. By definition of $\Cal R_M$ and $\Cal R_N$
we have $X_i=U_iX'_i$ and $Z_i=(U_i-V_i)X'_i$, $i=1,2$, where
$U_i\to V_i$ is an elementary $K$-reduction of types $(i)$--$(vi)$
if $Z_i\in\Cal R_N$ and $U_i=x_n^{\pm1}$, $V_i=u_\pm$, $X'_i\in F_n^+$ if $Z_i\in\Cal R_M$.

Following [\refFunar] and [\refBF], we represent
such sequences $W_1,\dots,W_m$ by diagrams. A {\it diagram} is
a union of mutually transversal curves
in the cylinder $S^1\times[0,1]$, each curve being
labeled by a letter $x_i^{\pm1}$.
In pictures we represent the cylinder by a rectangle whose vertical
sides are supposed to be identified, so, the fibers of the 
projection $\pr_2:S^1\times[0,1]\to[0,1]$
will be called {\it horizontal circles}.
Each curve is {\it monotone}, i.~e.,
its projection onto $[0,1]$ is bijective.
We say that a diagram is {\it admissible} if two curves 
labeled by $x_i^{\pm1}$ and $x_j^{\pm1}$ may cross only if $|i-j|\ge2$.
The words $W_i$ (up to cyclic permutation) are read on horizontal circles.

We say that curves $\Gamma_1,\dots,\Gamma_m$ form a
{\it bunch of parallel curves} or just a {\it bunch} if they
are pairwise disjoint and all the crossings
lying on $\bigcup\Gamma_i$ can be covered by disks whose intersections
with the diagram are as in Figure \figBunch\ up to symmetry.

In our case, the first and the last word of the sequence are $X_1$ and $X_2$.
So, on the boundary of the cylinder we indicate 
(by a bold line) segments corresponding to $U_1$ and $U_2$.
%
%
As in [\refFunar] and [\refBF], a diagram is called {\it interactive}
if it contains a curve which joins the bold segments.
We also say that a curve is {\it active} if it meets at least one
bold segment.


\proclaim{ Step 1 }
If all active curves form a single bunch all whose ends are on the bold
segments, then $Z_1-Z_2\in H$.
\endproclaim

In this case we have $U_1=U_2$.
Let $V_1=\rr(U_1)=\sum c_j W_j$, $c_j\in A$, $W_j\in F_\infty^+$.
For each $j$ we consider the diagram obtained from the initial diagram
by replacing the bunch of active curves by a bunch of curves labeled by
$W_j$. If a curve crosses the bunch, its label commutes with all letters
occurring in $U_1$, hence it commutes with all letters in $W_j$, i.~e.,
the new diagram is admissible and it defines a congruence
$W_jX_1'\equiv W_jX_2'\mod H$. Hence
(recall that $X_1-X_2\in H$) we have
$Z_1-Z_2=(X_1-V_1X_1')-(X_2-V_1X_2')
  \equiv 
  V_1X_2'-V_1X_1'=\sum c_jW_j(X_2'-X_1')\equiv0\mod H$.


\proclaim{ Step 2 }
If $Z_1,Z_2\in\Cal R_M$, then $Z_1-Z_2\in H$.
\endproclaim

In this case there is only one active curve, so we apply
the result of Step 1.


\proclaim{ Step 3 }
If the diagram is non-interactive, then $Z_1-Z_2\in H+R_{[w-1]}$.
\endproclaim

Due to Step 2, we may suppose that $Z_1\in\Cal R_N$.
Then $U_1=x_n^{\eps_1} x_{n-1}^{\eps_2} x_n^{\eps_3}$ with
$\eps_1,\eps_3\in\{-1,1\}$ and $\eps_2\in\{-1,0,1\}$.

\midinsert
\epsfxsize=90mm
\centerline{\epsfbox{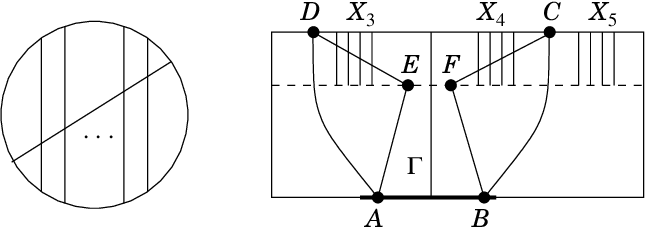}}
\botcaption{
              \hbox{Figure \figBunch
	      \hskip 35mm
              Figure \figNonInteract
	      \hskip 14mm}
}
\endcaption
\endinsert

Let $A$ and $B$ be the points on the lower bold segment
that correspond to the letters $x_n^{\eps_1}$ and $x_n^{\eps_3}$ of
$U_1$ and let $AD$ and $BC$ be the corresponding active curves
(see Figure \figNonInteract).
They cut the cylinder into two halves. Let $Q$ be that half whose side $AB$
is contained in the bold segment (the quadrangle $ABCD$ in Figure
\figNonInteract).

Let $\Gamma$ be the curve outcoming from $U_1$
and labeled by $x_{n-1}^{\eps_2}$ if $\eps_2\ne 0$ or a generic
monotone curve in $Q$ if $\eps_2=0$.
Let us choose a horizontal circle (the dashed line in Figure \figNonInteract)
so that all crossings are below it and let us choose points $E$ and $F$
on it so that the segment $EF$ which crosses $\Gamma$ has no other
intersections with the diagram.
We may suppose that the intersection of the diagram with the upper
half-cylinder (above $EF$) is a union of segments of vertical lines.

Let $\Delta$ be the diagram obtained by replacing $AD$ and $BC$ with
monotone curves $AED$ and $BFC$ where $ED$, $FC$ are straight
line segments and $AE$, $BF$ are curves in $Q$ which are chosen
so close to $\Gamma$ that the active curves outcoming from $U_1$
form a bunch in the lower half-cylinder (below $EF$).
The label of any curve $\Gamma'\ne\Gamma$ entering $Q$ is not $x_i^{\pm1}$
with $|n-i|\le 1$ (indeed, since $\Gamma'$ attains the
lower boundary outside the bold segment, it crosses $AD$ or $BC$).
Hence $\Delta$ is admissible.

Let $Y$ be the word read from $\Delta$ along the circle $EF$.
The bunch of active curves in the lower half-cylinder ensures that
$Y=U_1Y'$ and the result of Step 1 yields
$$
   Z_1\equiv (U_1 - \rr(U_1))Y'\mod H.         \eqno(\eqNonInteractOne)
$$

Now, let us study the upper part of $\Delta$ (above $EF$).
All the possible crossings in this part are
on $ED$ and $FC$. Hence, up to cyclic permutation, we have
$X_2= x_n^{\eps_1} X_3 x_{n-1}^{\eps_2} X_4 x_n^{\eps_3} X_5$
and $Y=U_1 Y' = U_1 X_4 X_5 X_3$
(see Figure \figNonInteract).
Since the diagram is not interactive, $U_2$ is a subword
of one of $X_3,X_4,X_5$, hence the active curves outcoming
from $U_2$ form a bunch and $Y'=Y_1U_2Y_2$, i.~e.,
$Y=U_1Y_1U_2Y_2$, $Y_1,Y_2\in F_\infty^+$. Hence,
by Step 1, we have
$$
  Z_2\equiv U_1Y_1(U_2-\rr(U_2))Y_2\mod H.    \eqno(\eqNonInteractTwo)
$$
We have also
$$
  U_1Y_1\rr(U_2)Y_2\equiv\rr(U_1)Y_1\rr(U_2)Y_2
  \equiv \rr(U_1)Y_1U_2Y_2 \mod R_{[w-1]}.
$$
Combining this with (\eqNonInteractOne) and (\eqNonInteractTwo), we obtain
$$
  Z_1\equiv (U_1 - \rr(U_1))Y_1U_2Y_2
     \equiv U_1Y_1(U_2 - \rr(U_2))Y_2\equiv Z_2\mod H+R_{[w-1]}.
$$
\medskip

\proclaim{ Step 4 }
Consider the open intervals obtained after removing all endpoints
of all active curves. If at least one of the words
corresponding to these intervals is not almost $K$-reduced
{\rm(}see the definition in \S\sectReduc\/{\rm)}, then
$Z_1-Z_2\in H + R_{[w-1]}$.
\endproclaim

Suppose that the word which is not almost $K$-reduced is a subword $Y$
of $X_2$. Since it is disjoint from the active curves, we can write
$X_2=U_2 X_3 Y X_4$. The fact that $Y$ is not almost $K$-reduced means
that there exists a sequence $Y=Y_0\to Y_1\to\dots\to Y'U_3Y''$ of exchanges of
commuting letters such that $U_3$ is the left hand side of
an elementary replacement of type
$(i)$--$(vi)$. The fact that $Y$ does not meet any active curve means that
the diagrams corresponding to the both chains
$$
\split
   X_1\to\dots\to &X_2=U_2X'Y_0X''\to U_2X'Y_1X''\to\dots\to U_2X'(Y'U_3Y'')X'',
\\
                  &X_2=U_2X'Y_0X''\to U_2X'Y_1X''\to\dots\to U_2X'(Y'U_3Y'')X''
\endsplit
$$
are non-interactive. By Step 3 this implies
$Z_1\equiv Z_3\equiv Z_2\mod H + R_{[w-1]}$
where $Z_3=U_2X'Y'(U_3-\rr(U_3))Y''X''$.

\proclaim{ Step 5 }
Suppose that $Z_1\in\Cal R_N$ and the diagram is interactive but not as in Step 1.
Then the active curves
are arranged up to symmetry either as in Figure \figInteract.1
or as in Figure \figInteract.2 where each of the dashed lines may or may
not be included into the diagram, $n\ge 1$.
\endproclaim

\midinsert
\epsfxsize=109mm
\centerline{
$x_n^{\pm1}\;x_{n-1}^{\pm1}\;x_n^{\pm1}\;x_{n-1}^{\pm1}\;x_n^{\pm1}
 \qquad\qquad\qquad\qquad
 x_{n-1}^{\pm1}\;x_n^{\pm1}\;x_{n-1}^{\pm1}\;x_n^{\pm1}$}
\centerline{\epsfbox{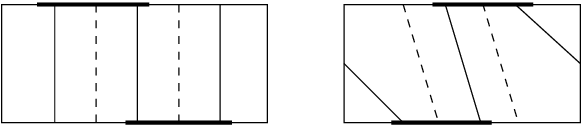}}
\centerline{
$x_n^{\pm1}\;x_{n-1}^{\pm1}\;x_n^{\pm1}\;x_{n-1}^{\pm1}\;x_n^{\pm1}
 \qquad\qquad\qquad\qquad
 x_n^{\pm1}\;x_{n-1}^{\pm1}\;x_n^{\pm1}\;x_{n-1}^{\pm1}$}
\botcaption{ Figure \figInteract.1 \hskip40mm Figure \figInteract.2  }
\endcaption
\endinsert

Indeed, we draw the curves adjacent to one of the bold segments
and we try all the possibilities to complete the picture to an
admissible diagram. 
%
%
It is easy to see that the pictures that could arise this way are the
two pictures from the statement of Step 5.

\medskip

\proclaim{ Step 6 }
If the active curves are as in Figure \figInteract.1, then
$Z_1-Z_2\in H+R_{[w-1]}+\tau(J_4)$.
\endproclaim

Suppose that the active curves are as in Figure \figInteract.1
(the bottom boundary corresponds to $X_1$). Then
$U_1=x_n^{\eps_4}x_{n-1}^{\eps_5}x_n^{\eps_6}$,
$U_2=x_n^{\eps_1}x_{n-1}^{\eps_2}x_n^{\eps_4}$,
$X_2=U_2 Y x_{n-1}^{\eps_5} X_3 x_n^{\eps_6} X_4$ where
$\eps_1,\eps_4,\eps_6=\pm1$ and $\eps_2,\eps_5\in\{-1,0,1\}$.

We begin as in Step 3.
Let $Q$ be the curvilinear quadrangle adjacent to the lower bold
segment and bounded by the active $x_n$-curves outcoming from $U_1$.
Let $C_1$ be a horizontal circle
such that the part of the diagram above $C_1$ is a union of segments
of vertical lines. Let $\Gamma$ be either the $(x_{n-1})$-curve
outcoming from $U_1$ (if it exists) or just a generic monotone curve in $Q$.
Then we push the $x_n$-curves inside the domain $Q$ from its boundary
so that they form (together with $\Gamma$) a bunch below $C_1$,
and so that the portions of the pushed curves above $C_1$
are segments of straight lines (see Figure \figInteractA.1).

\midinsert
\centerline{
   $\;\;\;U_2\;\;\; Y  \;\; x_{n-1}^{\eps_5} \; X_3 \;\; x_n^{\eps_6} \;\; X_4$
     \hskip15mm
    $\;U_2\;\;\; Y_1\;\; x_{n-2}^{\eps_3} Y_2  \;\, x_{n-1}^{\eps_5}
        X_3 \;\; x_n^{\eps_6}\;\, X_4$
}
\epsfxsize=110mm
\centerline{\epsfbox{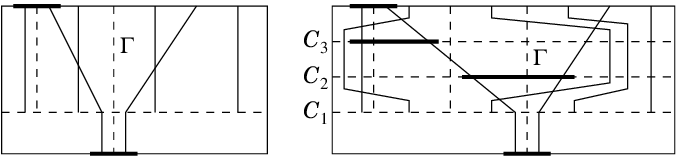}}
\centerline{$U_1$ \hskip62mm $U_1$ \hskip 5mm}
\botcaption{
    Figure \figInteractA.1
    \hskip 40mm
    Figure \figInteractA.2
    \hskip 15mm
 }
 \endcaption
\endinsert

Since all curves outcoming from $Y$ cross an $x_n$-curve,
$Y$ does not contain $x_i^{\pm1}$ for $n-1\le i\le n+1$.
By Step 4, we may suppose that $Y$ is almost $K$-reduced, hence
$Y$ has at most one occurrence of $x_{n-2}^{\pm1}$,
i.~e., $Y=Y_1 x_{n-2}^{\eps_3}Y_2$ with $\eps_3\in\{-1,0,1\}$ and
$Y_1,Y_2$ do not contain $x_i^{\pm1}$ for $n-2\le i\le n+1$.

We choose horizontal circles $C_2$ and $C_3$ so that
the intersection point of the $x_n$-curve and the $(x_{n-2})$-curve (if it exists)
is between them and
we modify the diagram as it is shown in Figure \figInteractA.2.
If we apply the result of Step 1 to the part of the diagram which
is below $C_2$ and to that which is above $C_3$, we obtain:
$$
\xalignat2
   &Z_1\equiv Y_1x_n^{\eps_1} x_{n-1}^{\eps_2} x_{n-2}^{\eps_3}
   \big(U_1-\rr(U_1)\big)Y_2X_3X_4 \mod H
  &&\text{(below $C_2$),}
\\
   &Z_2\equiv Y_1\big(U_2-\rr(U_2)\big) x_{n-2}^{\eps_3}x_{n-1}^{\eps_5}
   x_n^{\eps_6}Y_2X_3X_4 \mod H.
  &&\text{(above $C_3$).}
\endxalignat
$$
Hence $Z_1-Z_2\equiv X'\sh^{n-3}Y'\mod H$
where $X'=Y_2X_3X_4Y_1$ and
$$
  Y'=
    \rr(x_3^{\eps_1}x_2^{\eps_2}x_3^{\eps_4})\,x_1^{\eps_3}
        x_2^{\eps_5}x_3^{\eps_6} -
    x_3^{\eps_1}x_2^{\eps_2}x_1^{\eps_3}\,
    \rr(x_3^{\eps_4}x_2^{\eps_5}x_3^{\eps_6}).
$$
If $\eps_2\ne0$, then $\rr(Y')\in J_4$ by Condition (J1) of the
definition of $J_4$. Thus, using Lemma \lemDiffiOne\ and observing that
$\wt(X'\sh^{n-3}Y')<w$,
we obtain
$$
  Z_1-Z_2\equiv X'\sh^{n-3}Y'
       \equiv X'\sh^{n-3} \rr(Y')
       \equiv 0 \mod H+R_{[w-1]}+\tau(J_4).
$$
If $\eps_2=0$, then
$X'\sh^{n-3}Y'\equiv X'\sh^{n-3}(x_1^{\eps_3}Y'')\mod H$
where $\rr(Y'')\in J_4$, thus
$$
  Z_1-Z_2\equiv X'\sh^{n-3}(x_1^{\eps_3}Y'')
       \equiv X'x_{n-2}^{\eps_3}\sh^{n-3} \rr(Y'')
       \equiv 0 \mod H+R_{[w-1]}+\tau(J_4).
$$

\proclaim{ Step 7 }
If the active curves are as in Figure \figInteract.2, then
$Z_1-Z_2\in H+R_{[w-1]}+\tau_N(L)$.
\endproclaim

\midinsert
\centerline{ $X_3\;\;x_{n-1}^{\eps_2}\;\; X_4  \;\;\; U_2$  }
\vskip-1pt
\epsfxsize=38mm
\centerline{\epsfbox{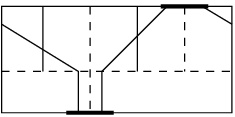}}
\vskip-2.5pt
\centerline{ $U_1\qquad\;$ }
\botcaption{
    Figure \figInteractB
 }
\endcaption
\endinsert

Again, as in the beginning of Steps 3 and 6, we transform the
diagram as in Figure \figInteractB\ and we obtain
$$
  Z_1\equiv X_3(U_1-\rr(U_1))X_4 x_{n-1}^{\eps_4}
  \;\;\text{and}\;\;
  Z_2\equiv X_3x_{n-1}^{\eps_2}X_4(U_2-\rr(U_2))\mod H
$$
where
$U_1=x_n^{\eps_1}x_{n-1}^{\eps_2}x_n^{\eps_3}$,
$U_2=x_n^{\eps_3}x_{n-1}^{\eps_4}x_n^{\eps_1}$,
$\eps_1,\eps_3=\pm1$, $\eps_2,\eps_4\in\{-1,0,1\}$.
Hence
$$
  Z_1-Z_2\equiv
    X_3x_{n-1}^{\eps_2}X_4\rr(U_2) - X_3\rr(U_1)X_4 x_{n-1}^{\eps_4}
    \mod H
                                                      \eqno(\eqInterB)
$$
Note that $x_i^{\pm1}$ for $n-1\le i\le n+1$ does not occur
in $X_j$, $j=3,4$. Indeed, if it does, then the diagram
$x_i^{\pm1}$-curve starting at a point from the upper circle corresponding to the letter $x_i^{\pm1}$
in $X_j$
cannot attain the opposite side of the cylinder
outside the bold segment because it cannot cross the $x_n$-curves.
Thus,
$$
  X_3x_{n-1}^{\eps_2}X_4\rr(U_2)\equiv
  X'_3x_{n-1}^{\eps_2}X'_4\rr(U_2)X_5
  \equiv X'_3x_{n-1}^{\eps_2}X'_4\rr(U_2) \tau(X_5)
  \mod H+R_{[w-1]}
$$
where $X'_3,X'_4\in F_{n-1}^+$ and $X_5\in\sh^{n+1} F_\infty^+$
(the same for other term in (\eqInterB)). So, replacing, if necessary,
$X_j$ by $X'_j$ ($j=3,4$), we may assume that $X_3,X_4\in F_{n-1}^+$.
Then we can pass from (\eqInterB) to
$$
\split
  Z_1-Z_2&\equiv
    X_3x_{n-1}^{\eps_2}X_4\tau_{n+1}(U_2) - X_3\tau_{n+1}(U_1)X_4 x_{n-1}^{\eps_4}
    \mod H + R_{[w-1]}
\\&
    =\varphi_{n-1,X_3,X_4}(Y)
\endsplit
$$
where
$Y=x_1^{\eps_2}\otimes\tau_3(x_2^{\eps_3}x_1^{\eps_4}x_2^{\eps_1})
-\tau_3(x_2^{\eps_1}x_1^{\eps_2}x_2^{\eps_3})\otimes x_1^{\eps_4}
\in L$.
Thus  if $n>1$, then the result follows from Lemma \lemDiffiTwo.
If $n=1$, then $\eps_2=\eps_4=0$, $X_3=X_4=1$, and (\eqInterB)
yields $Z_1-Z_2\equiv 0$.
%
\qed
\enddemo


\proclaim{ Lemma \lemDiffiThree}
$R_{[w]}\cap AF^+_{[w-1]}\subset H+R_{[w-1]}+\tau(J_4)+\tau_N(L)$
\endproclaim

\demo{ Proof }
For $Z\in R_{[w]}\cap AF^+_{[w-1]}$, let $m=m(Z)$ be the minimal number
such that $Z\equiv c_1Z_1+\dots+c_mZ_m\mod H+R_{[w-1]}$ with $c_i\in A$,
$Z_i\in\Cal R_N\cup\Cal R_M$, $\wt Z_i=w$.
To prove that $Z\in H+R_{[w-1]}+\tau(J_4)+\tau_N(L)$, we use
induction on $m$. The statement is trivial for $m=0$.

Suppose that $m>0$ and the statement is true for all smaller values of $m$.
Let $X_i$ be the leading monomial of $Z_i$, i.~e.,
$X_i\in F_\infty^+$ and $\wt(Z_i-X_i)<w$. Then $\sum c_i X_i\equiv 0\mod H$.
The term $c_mX_m$ of this congruence should be cancelled by other terms.
Hence there exists $j<m$ such that $X_m-X_j\in H$.
Then $c_m(Z_m-Z_j)\in H+R_{[w-1]}+\tau(J_4)+\tau_N(L)$ by Lemma \lemPent\
and $Z-c_m(Z_m-Z_j)\in H+R_{[w-1]}+\tau(J_4)+\tau_N(L)$ by the induction
hypothesis.
\qed\enddemo

\proclaim{ Proposition \propDiffi } $I\subset \tau(J_4)+\tau_N(L)$.
\endproclaim

\demo{ Proof }
Let $I'=\tau(J_4)+\tau_N(L)$.
Since $I=R\cap A$ and $R=\bigcup_w R_{[w]}$, it is enough to prove
that $R_{[w]}\cap A\subset I'$ for any $w$.
For $w=0$ we have $R_{[0]}=0$, hence $R_{[0]}\cap A\subset I'$.
Suppose that $R_{[w-1]}\cap A\subset I'$.
Let $Z\in R_{[w]}\cap A$. 
Since $R_{[w]}\cap A\subset R_{[w]}\cap AF^+_{[w-1]}$,
by Lemma \lemDiffiThree\ we have
$Z\in H+R_{[w-1]}+I'$, i.~e., $Z=Z_H+Z'+Z_0$ with $Z_H\in H$, $Z'\in R_{[w-1]}$,
$Z_0\in I'$. Since $Z_H=Z-Z'-Z_0\in AF^+_{[w-1]}$ and $H$ is homogeneous, we
have $Z_H\in R_{[w-1]}$, thus $Z=Z''+Z_0$ with $Z''=Z'+Z_H\in R_{[w-1]}$
and $Z_0\in I'$. Since $Z''=Z-Z_0\in A$, we have
$Z''\in R_{[w-1]}\cap A$ and by the induction hypothesis we
obtain $Z''\in I'$ whence $Z=Z''+Z_0\in I'$. Thus $R_{[w]}\cap A\subset I'$.
\qed
\enddemo

Our main Theorem follows from Propositions \propEasy\
and \propDiffi.

\head 4. The link invariant $P_{\alpha,\beta;k}$
\endhead

\subhead\sectInv. Some general properties
\endsubhead
As we mentioned in Remark \remInv, the Markov trace on $AK_\infty$ defines an invariant of oriented links
$P=P_{\alpha,\beta;k}$ which takes its values in $\tilde A/I$ where
$\tilde A=k[u^{\pm1/2},v^{\pm1/2}]$.
If $u$ and $v$ are not zero divisors (this is so in all the cases computed so far;
see Remark \remTransv), then $A$ embeds into $\tilde A$.

\proclaim{ Proposition \propInv }
a). $P(\text{\rm unknot})=1$.

\smallskip\noindent
b). $P$ satisfies the skein relation coming from (\eqRelOne):
$$
   u^{3/2} P(\,\epsfysize=8pt\lower1pt\hbox{{\epsfbox{skein2.eps}}}\,)
 - \alpha u v^{1/2} P(\,\epsfxsize=8pt\lower1pt\hbox{{\epsfbox{skein1.eps}}}\,)
 + \beta u^{1/2}vP(\,\epsfxsize=8pt\lower1pt\hbox{{\epsfbox{skein0.eps}}}\,)
 - v^{3/2}P(\,\epsfxsize=8pt\lower1pt\hbox{{\epsfbox{skein-.eps}}}\,)           \eqno(\eqSkein)
$$
Similarly, (\eqRelTwo) and (\eqRelTwoBar)
yield skein relations for tangles with $6$ endpoints.

\smallskip\noindent
c). $P(L_1\sqcup L_2) = (uv)^{-1/2} P(L_1)P(L_2)$ (disjoint union).

\smallskip\noindent
d). $P(L_1\# L_2) = P(L_1)P(L_2)$ (connected sum).

\smallskip\noindent
e). If $\bar L$ is the mirror image of $L$, then $P_{\alpha,\beta}(\bar L)(u,v)=P_{\beta,\alpha}(L)(v,u)$.

\smallskip\noindent
f). If $L_2$ is obtained from $L_1$ by a mutation (i.~e., a rotation of a tangle with $4$ endpoints
by $180^\circ$ such  that no endpoint is fixed),
then $P(L_1)=P(L_2)$.

\smallskip\noindent
g). $P(L)\equiv 1\mod \big(u-1,\,v-1,\,\alpha-\beta,\,2(\alpha-2)^2\big)$.

\smallskip\noindent
h). If $L$ can be represented by an $n$-braid, then $P(L)$ can be represented by a Laurent polynomial
$f\in\tilde A$ such that $1-n\le\deg_{u,v}(T)\le 0$ for any monomial $T$ of $f$.
\endproclaim

\demo{ Proof } a) -- e) and h). Immediate from the definition of $P$.

\smallskip
f). Mutant links $L_1$ and $L_2$ can be represented by braids $Z_1=X\sh^{n-2}Y$ and
$Z_2=X\sh^{n-2}(\sigma_1 Y \sigma_1^{-1})$ 
respectively where $X\in B_n$, $Y\in B_m$.
Hence we have
$t(Z_1) = t(\tau_{n+1}\circ\tau_{n+2}\circ\dots\circ\tau_{n+m-2}(Z_1)) = t(X\sh^{n-2}Y_2)$
where $Y_2=\tau_3\circ\dots\circ\tau_m(Y)\in AB_2$ and, similarly,
$t(Z_2) = t(X\sh^{n-2}(\sigma_1 Y_2 \sigma_1^{-1}))$.
It remains to note that
$Y_2=\sigma_1 Y_2 \sigma_1^{-1}$
because the group $B_2$
is abelian.

\smallskip
g). If we set $x=y=\bar x=\bar y=1$ in (\eqRelOne) and (\eqRelTwo), then
we obtain identities modulo the ideal $(\alpha-\beta,\,2(\alpha-2)^2)$.
\qed
\enddemo

From now on we assume that $k=\bold k[\alpha,\beta]$ where $\bold k$ is a commutative ring
and each of $\alpha$ and $\beta$ is either zero or transcendent over $\bold k$.
For a monomial in $\alpha,\beta,u,v$ we define its degree modulo 3 (denoted by $\deg_3$) by
setting $\deg_3(u)=\deg_3(\alpha)=1$ and $\deg_3(v)=\deg_3(\beta)=-1$.
We denote the $\deg_3$-homogeneous component of $A$ (resp. of $\tilde A$) of degree $d$
by $A_d$ (resp. by $\tilde A_d$). So, we have $A=A_0\oplus A_1\oplus A_2$ and
$\tilde A=\tilde A_0\oplus \tilde A_{1/2}\oplus\dots\oplus\tilde A_{5/2}$.
The following fact follows immediately from the definitions.

\proclaim{ Proposition \propPeriod }
The ideal $I$ is $\deg_3$-homogeneous. $P(L)\in\tilde A_0/I$ for any $L$.
\qed
\endproclaim


\subhead \sectNF. Normal form of elements of $k[u^{\pm1},v^{\pm1}]/I$
\endsubhead
In fact, the square roots of $u$ and $v$ in the definition of the invariant $P$
are needed only for writing the skein relation (\eqSkein) in a nice form.
Otherwise, for a $\mu$-component link $L$ represented by an $n$-braid $X$, we can set
$$
   P(L)=u^{(p-n-e)/2}v^{(p-n+e)/2}t(X) \quad\text{where}\quad
   p=\cases 1, &\text{ $\mu$ is odd,}\\ 0, &\text{ $\mu$ is even}\endcases        \eqno(\eqDefP)
$$
which ensures that $P(L)$ belongs to $k[u^{\pm1},v^{\pm1}]/I$ for any link $L$.
Note that Proposition \propPeriod\ still holds for $P(L)$ defined by (\eqDefP).
So, from now on we forget about fractional powers of $u$ and $v$ and we discuss the normal
form of elements of the ring $k[u^{\pm1},v^{\pm1}]$. 

Assume that $k=\bold k[\alpha,\beta]$ as in the previous subsection.
We have a natural identification of $k[u^{\pm1},v^{\pm1}]$ with $\bar A/\bar I$ where
$\bar A=A[\bar u,\bar v]$ ($\bar u$ and $\bar v$ are new independents variables) and
$\bar I=I+(u\bar u-1,v\bar v-1)$. So, to define a normal form in $k[u^{\pm1},v^{\pm1}]$,
it is enough to compute a Gr\"obner base $\bar{\Cal G}$ of $\bar I$.

We assume that the monomial order is chosen so that $\bar u$ and $\bar v$ are greater than any
monomial in $\alpha,\beta,u,v$. In this case, the following conditions are equivalent:
\roster
\item $A/I$ embeds into $\bar A/\bar I\cong k[u^{\pm1},v^{\pm1}]/I$;
\item $u$ and $v$ are not zero divisors in $A/I$;
\item $\bar{\Cal G}$ contains a Gr\"obner base of $I$.
\endroster
Condition (3) holds in all the cases computed so far (see Remark \remTransv).

If $\bold k$ is a field, then the normal form in $\bar A/\bar I$ defined by $\bar{\Cal G}$ is evident:
it is just a $\bold k$-linear combination of monomials
which are not divisible by the leading terms of elements of $\bar{\Cal G}$.
In the case when $\bold k$ is $\Bbb Z$ or $\Bbb Z/m\Bbb Z$, the situation is more delicate.
For any monomial $T$, we have to fix canonical representatives in $\bold k$ for elements
of $\bold k/I(T)$ where $I(T)$ is the ideal of $\bold k$ consisting of the leading coefficients
of elements of $\bar A$ whose leading monomial is $T$. When $I(T)\ne 0$, we choose them
in $\{0,\dots,m(T)-1\}$ where $m(T)$ is the positive generator of $I(T)$. Thus the
normal form is a $\bold k$-linear combination of monomials $T$ with coefficients
belonging to $\{0,\dots,m(T)-1\}$. Note that $m(T)$ is the gcd of the leading coefficients of
those elements $\bar{\Cal G}$ whose leading monomials divide $T$.



\subhead \sectGBQ. The case $\bold k=\Bbb Q$ or $\Bbb F_p$ and $\beta=0$
\endsubhead

\proclaim{ Proposition \propGBQ } Let the notation be as in the respective parts of Corollary \corBetaZero\
and let $\bar{\Cal G}$ be the reduced Gr\"obner base of $\bar I$ with respect to the lexicographic
order such that $\bar v>\bar u>v>u>\alpha$. Then:

\smallskip
a). {\rm($\bold k=\Bbb F_2$)}. $\bar{\Cal G}=\Cal G\cup\{u\bar u+1,\,v\bar v+1,\,
     \alpha^2(\bar u+ u^2+\alpha u),\,   
     \alpha^2(\bar v+ u+\alpha),\break       
     u^3\bar v+v^2+u+\alpha(uv+1+\alpha u^2+\alpha^2u),\, 
     \bar u\bar v+\bar uv^2+1+\alpha(\bar u+ v+\alpha u+\alpha^2)\}$.
\if01{
I[5]=Ua2+v2ua2+vu2a3+ua3
I[6]=Uu+1
I[7]=Va2+U2a2
I[8]=Vu3+Vua2+V+v2+vua+u+a
I[9]=Vv+1
I[10]=VU+Vu2+Va2+Uv2+Uvua+Uu+Ua
5: u2a2+ua3
7: ua2+a3
8: V+v2+vua+u2a2+ua3+u+a
10: Vu2+Uv2+Ua+va+ua2+a3+1
}\fi

\smallskip
b). {\rm($\bold k=\Bbb F_3$)}. $\bar{\Cal G}=\Cal G\cup\{\bar u+\alpha^2 u^3-u^2+\alpha u+\alpha^2,\,
                              \bar v+\alpha u^3+\alpha^2 u^2-u\}$.

c). {\rm($\bold k=\Bbb Q$ or $\Bbb F_p$, $p\not\in\{2,3,37\}$)}. $\bar{\Cal G}=\Cal G\cup\{f_6,f_7\}$ where

\medskip
$
  \;f_6=\bar u + {7\over12}\alpha^2 u^3 + ({7\over8}\alpha^3+6) u^2 - ({9\over8}\alpha^3+{20\over3})\alpha u
  + {1\over48}\alpha^8 + {163\over288}\alpha^5+{28\over9}\alpha^2,
$

\medskip
$
  \;f_7=\bar v + {35\over6}\alpha u^3-({69\over32}\alpha^3+{71\over4})\alpha^2 u^2
-({75\over8}\alpha^3+{152\over3})u+{155\over96}\alpha^7+{2867\over144}\alpha^4+{433\over9}\alpha.
$

\medskip
d).  {\rm($\bold k=\Bbb F_{37}$)}. $\bar{\Cal G}=\Cal G\cup\{f_{6,37},f_{7,37}\}$ where
$$
\split
  &f_{6,37}=\bar u-3\alpha u^4-6\alpha^2u^3+(6\alpha^3-2)u^2-10(\alpha^3-1)\alpha u-2\alpha^8+15\alpha^5+14\alpha^2,\\
  &f_{7,37}=\bar v+14u^4+16\alpha u^3-(12\alpha^3u^2+11)\alpha^2u^2+(12\alpha^3-1)u-12\alpha^7-6\alpha^4-\alpha.
\endsplit
$$

\if01{
\medskip
e).
{\def\U{\bar u}\def\a{\alpha}\def\au{\alpha u}\def\V{\bar v}
$\Cal G\cup\{u\U-1,\,v\V-1,\,
4\a\U-4\au^5+3u^2\a^4+8\au^2+8\a^2u-2\a^3,\,     
\a^3\U+8\U+11\a^3u^2+8u^2+8\a^2,\,               
8\U v+6\a^2\U-2\a^2v+4\au^3+8\a^2u^2+8u+4\a,\,   
4\a\V+4\au^2v+8\au,\,                            
\a^3\V+8\V+8\a^2u^2+9\a^3u+8u,\,                 
2\a^2u\V+8\V+4\au v+2\a^2v-4\au^3+8\a^2u^2+8u,\, 
4u^2\V+2\au\V+7\a^2\V+2\a v+7\a^2u+9\a^3+4,\,    
2\au^2\V-\a^2u\V-4\V-4\au v+7\a^2v-3\a^3u+4u+2\a,\,         
\a^2u^2\V+4u\V-2\a\V+4v-2\au+3a^2,\,                        
u^3\V-3\V+5v^2+3\au v-6\a^2v-7\a^2u^2+5\a^3u+5u+9\a^4+\a,\, 
\U\V+u^2\V-2\au\V+5\a^2\V+\U v^-7\a\U+\a v+8\au^2+4\a^2u-3  
\}$} is a Gr\"obner base of $\bar{I}$.
}\fi 
\endproclaim



Thus, in the setting of Corollary \corBetaZero(c) ($\bold k=\Bbb Q$ or $\Bbb F_p$, $p\not\in\{2,3,37\}$),
the normal form of elements of $\bar A/\bar I$ always belongs to $A$ and
we have $\bar A/\bar I\cong A/I$, in particular, $\dim_{\bold k}\bar A/\bar I=24$ (or $21$ if $\bold k=\Bbb F_7$).
By Proposition \propPeriod, the invariant $P(L)$ takes its values in $\tilde A_0/I$.
So, its normal form is a linear combination of the eight monomials indicated in the header line of Table 1
(without $u^3$ in the case $\bold k=\Bbb F_7$). The values of $P_{0,0,\bold k[\alpha]}$ for knots up to 9 crossings
are presented in Table 2 (the choice between the ``right'' knots $3_1,5_1,5_2,\dots$ and their mirror images
$\bar 3_1,\bar 5_1,\bar 5_2\dots$ is done according to the database ``The Knot Atlas'' {\tt http://katlas.org}).

The most interesting case is $\bold k=\Bbb F_p$ for $p=37$. Up to now this is the only case when
the invariant $P$ distinguishes knots with equal HOMFLY polynomials.
In this case the normal form has one more monomial: $\alpha^2 u^4$.

\topinsert
\noindent Table 1. $P_{\alpha,0;\Bbb Q}(K)$ for knots $K$ up to 9 crossings ($\bar K$ is the mirror of $K$)
\smallskip
\setbox1=\vbox{\offinterlineskip
\hrule
\def\o{&\omit&}
\def\tabA{height3.1pt
& \omit &
& \omit \o \omit \o \omit \o \omit \o \omit \o \omit \o \omit \o \omit &
& \omit \o \omit \o \omit \o \omit \o \omit \o \omit \o \omit \o \omit & \cr}
\def\NL{&\cr\tabA&\,}
\def\SP{&\omit&}
\def\VL{&&\,}
\def\m{\text{-}}
\def\same{\o\o\o\o\o\o\o}
\halign{&\vrule#&\strut#\hfill\cr
\tabA
& \lower7pt\hbox{\,$K$}&
&\multispan{15}\hfill coefficients of $P_{\alpha,0;\Bbb Q}( K)$ \hfill&
&\multispan{15}\hfill coefficients of $P_{\alpha,0;\Bbb Q}(-K)$ if differ\hfill&\cr
& \omit &
&\ $1$\o $\alpha^3$ \o $\alpha^6$ \o $\alpha^2u$ \o $\alpha^5u$ \o $\alpha u^2$ \o $\!\alpha^4\hskip-0.7pt u^2$ \o $u^3$ &
&\ $1$\o $\alpha^3$ \o $\alpha^6$ \o $\alpha^2u$ \o $\alpha^5u$ \o $\alpha u^2$ \o $\!\alpha^4\hskip-0.7pt u^2$ \o $u^3$ &\cr
\tabA
\noalign{\hrule}
\tabA
&\, $0_1$&&\ 1 \o 0 \o 0 \o 0 \o 0 \o 0 \o 0 \o 0 &&\same
\NL $3_1$\VL
${\m20\over 9}$\SP ${\m61\over 36}$\SP ${\m1\over 24}$\SP ${\m44\over 3}$\SP
${\m5\over 2}$\SP ${26\over 1}$\SP ${27\over 8}$\SP ${\m35\over 3}$\VL
${\m145\over 81}$\SP ${\m899\over 324}$\SP ${\m317\over 864}$\SP
${605\over 108}$\SP ${275\over 288}$\SP ${\m20\over 3}$\SP ${\m23\over 24}$\SP ${56\over 27}$
\NL $4_1$\VL
${\m10\over 9}$\SP ${25\over 9}$\SP ${43\over 96}$\SP ${\m127\over 12}$\SP
${\m49\over 32}$\SP ${10\over 1}$\SP ${5\over 4}$\SP ${\m7\over 3}$\VL
\same
\NL $5_1$\VL
${22\over 3}$\SP ${965\over 12}$\SP ${365\over 32}$\SP ${\m1221\over 4}$\SP
${\m1449\over 32}$\SP ${367\over 1}$\SP ${375\over 8}$\SP ${\m119\over 1}$\VL
${\m796\over 729}$\SP ${\m3229\over 1458}$\SP ${\m1565\over 7776}$\SP ${4769\over 972}$\SP
${659\over 2592}$\SP ${\m55\over 9}$\SP ${\m1\over 72}$\SP ${791\over 243}$
\NL $5_2$\VL
${34\over 9}$\SP ${311\over 9}$\SP ${491\over 96}$\SP ${\m1823\over 12}$\SP
${\m737\over 32}$\SP ${194\over 1}$\SP ${99\over 4}$\SP ${\m203\over 3}$\VL
${146\over 243}$\SP ${\m446\over 243}$\SP ${\m923\over 2592}$\SP ${2399\over 324}$\SP
${1313\over 864}$\SP ${\m82\over 9}$\SP ${\m59\over 36}$\SP ${203\over 81}$
\NL $6_1$\VL
${19\over 9}$\SP ${161\over 36}$\SP ${47\over 96}$\SP ${49\over 12}$\SP
${31\over 32}$\SP ${\m16\over 1}$\SP ${\m17\over 8}$\SP ${28\over 3}$\VL
${136\over 81}$\SP ${1799\over 324}$\SP ${22\over 27}$\SP ${\m437\over 27}$\SP
${\m179\over 72}$\SP ${50\over 3}$\SP ${53\over 24}$\SP ${\m119\over 27}$
\NL $6_2$\VL
${\m55\over 27}$\SP ${\m481\over 54}$\SP ${\m95\over 72}$\SP ${365\over 9}$\SP
${19\over 3}$\SP ${\m55\over 1}$\SP ${\m57\over 8}$\SP ${182\over 9}$\VL
${\m445\over 243}$\SP ${\m5\over 972}$\SP ${7\over 648}$\SP ${\m340\over 81}$\SP
${\m29\over 108}$\SP ${13\over 3}$\SP ${1\over 24}$\SP ${\m154\over 81}$
\NL $6_3$\VL
${\m157\over 81}$\SP ${\m1265\over 324}$\SP ${\m23\over 54}$\SP
${83\over 27}$\SP ${7\over 36}$\SP\ $0$\SP ${1\over 4}$\SP ${\m28\over 27}$\VL
\same
\NL $7_1$\VL
${178\over 3}$\SP ${2753\over 12}$\SP ${785\over 32}$\SP ${\m1257\over 4}$\SP
${\m1053\over 32}$\SP ${131\over 1}$\SP ${171\over 8}$\SP ${77\over 1}$\VL
${\m2936\over 2187}$\SP ${\m20393\over 4374}$\SP ${\m26419\over 23328}$\SP ${52159\over 2916}$\SP
${39997\over 7776}$\SP ${\m1819\over 81}$\SP ${\m3925\over 648}$\SP ${2569\over 729}$
\NL $7_2$\VL
${169\over 9}$\SP ${1265\over 9}$\SP ${449\over 24}$\SP ${\m1289\over 3}$\SP
${\m491\over 8}$\SP ${468\over 1}$\SP ${60\over 1}$\SP ${\m392\over 3}$\VL
${\m665\over 243}$\SP ${\m970\over 243}$\SP ${\m59\over 81}$\SP ${832\over 81}$\SP
${145\over 54}$\SP ${\m344\over 27}$\SP ${\m167\over 54}$\SP ${196\over 81}$
\NL $7_3$\VL
${1351\over 2187}$\SP ${2123\over 8748}$\SP ${2471\over 5832}$\SP ${\m2084\over 729}$\SP
${\m2461\over 972}$\SP ${89\over 27}$\SP ${737\over 216}$\SP ${2338\over 729}$\VL
${43\over 1}$\SP ${521\over 2}$\SP ${33\over 1}$\SP ${\m678\over 1}$\SP
${\m747\over 8}$\SP ${673\over 1}$\SP ${699\over 8}$\SP ${\m154\over 1}$
\NL $7_4$\VL
${\m28\over 81}$\SP ${\m989\over 324}$\SP ${\m155\over 216}$\SP ${326\over 27}$\SP
${13\over 4}$\SP ${\m410\over 27}$\SP ${\m815\over 216}$\SP ${77\over 27}$\VL
${223\over 9}$\SP ${6365\over 36}$\SP ${2291\over 96}$\SP ${\m6803\over 12}$\SP
${\m2621\over 32}$\SP ${636\over 1}$\SP ${651\over 8}$\SP ${\m560\over 3}$
\NL $7_5$\VL
${28\over 1}$\SP ${309\over 2}$\SP ${621\over 32}$\SP ${\m1601\over 4}$\SP
${\m1761\over 32}$\SP ${399\over 1}$\SP ${417\over 8}$\SP ${\m91\over 1}$\VL
${8650\over 2187}$\SP ${20987\over 8748}$\SP ${18569\over 23328}$\SP ${\m16697\over 2916}$\SP
${\m28751\over 7776}$\SP ${187\over 27}$\SP ${1051\over 216}$\SP ${2401\over 729}$
\NL $7_6$\VL
${20\over 27}$\SP ${1385\over 54}$\SP ${1093\over 288}$\SP ${\m4009\over 36}$\SP
${\m1603\over 96}$\SP ${139\over 1}$\SP ${141\over 8}$\SP ${\m427\over 9}$\VL
${\m542\over 243}$\SP ${\m1789\over 972}$\SP ${\m895\over 2592}$\SP ${1039\over 324}$\SP
${1081\over 864}$\SP ${\m43\over 9}$\SP ${\m115\over 72}$\SP ${49\over 81}$
\NL $7_7$\VL
${\m34\over 27}$\SP ${89\over 54}$\SP ${7\over 18}$\SP ${\m118\over 9}$\SP
${\m55\over 24}$\SP ${50\over 3}$\SP ${59\over 24}$\SP ${\m49\over 9}$\VL
${\m67\over 81}$\SP ${46\over 81}$\SP ${55\over 864}$\SP ${773\over 108}$\SP
${335\over 288}$\SP ${\m16\over 1}$\SP ${\m15\over 8}$\SP ${224\over 27}$
\NL $8_1$\VL
${\m35\over 9}$\SP ${\m286\over 9}$\SP ${\m14\over 3}$\SP ${424\over 3}$\SP
${43\over 2}$\SP ${\m184\over 1}$\SP ${\m47\over 2}$\SP ${196\over 3}$\VL
${\m173\over 243}$\SP ${1121\over 243}$\SP ${521\over 648}$\SP ${\m1457\over 81}$\SP
${\m659\over 216}$\SP ${172\over 9}$\SP ${26\over 9}$\SP ${\m392\over 81}$
\NL $8_2$\VL
${\m94\over 9}$\SP ${\m2165\over 36}$\SP ${\m701\over 96}$\SP ${1493\over 12}$\SP
${515\over 32}$\SP ${\m101\over 1}$\SP ${\m105\over 8}$\SP ${35\over 3}$\VL
${1928\over 729}$\SP ${2591\over 1458}$\SP ${4099\over 7776}$\SP ${\m5215\over 972}$\SP
${\m6301\over 2592}$\SP ${179\over 27}$\SP ${665\over 216}$\SP ${203\over 243}$
\NL $8_3$\VL
${226\over 81}$\SP ${812\over 81}$\SP ${1127\over 864}$\SP ${\m1307\over 108}$\SP
${\m437\over 288}$\SP ${2\over 3}$\SP ${1\over 12}$\SP ${133\over 27}$\VL
\same
\NL $8_4$\VL
${\m98\over 81}$\SP ${\m721\over 324}$\SP ${\m385\over 864}$\SP ${2905\over 108}$\SP
${1327\over 288}$\SP ${\m45\over 1}$\SP ${\m49\over 8}$\SP ${511\over 27}$\VL
${\m244\over 243}$\SP ${3245\over 486}$\SP ${2293\over 2592}$\SP ${\m5785\over 324}$\SP
${\m1723\over 864}$\SP ${43\over 3}$\SP ${25\over 24}$\SP ${\m259\over 81}$
\NL $8_5$\VL
${185\over 729}$\SP ${2443\over 2916}$\SP ${4015\over 7776}$\SP ${\m6967\over 972}$\SP
${\m7765\over 2592}$\SP ${245\over 27}$\SP ${203\over 54}$\SP ${98\over 243}$\VL
${\m148\over 9}$\SP ${\m1735\over 18}$\SP ${\m299\over 24}$\SP ${785\over 3}$\SP
${293\over 8}$\SP ${\m269\over 1}$\SP ${\m69\over 2}$\SP ${203\over 3}$
\NL $8_6$\VL
${\m151\over 27}$\SP ${\m5915\over 108}$\SP ${\m137\over 18}$\SP ${1745\over 9}$\SP
${343\over 12}$\SP ${\m228\over 1}$\SP ${\m117\over 4}$\SP ${644\over 9}$\VL
${\m101\over 729}$\SP ${1091\over 2916}$\SP ${\m35\over 243}$\SP ${\m413\over 243}$\SP
${323\over 324}$\SP ${4\over 3}$\SP ${\m19\over 12}$\SP ${\m644\over 243}$
\NL $8_7$\VL
${\m29\over 243}$\SP ${\m4489\over 972}$\SP ${\m257\over 324}$\SP ${1117\over 81}$\SP
${601\over 216}$\SP ${\m133\over 9}$\SP ${\m205\over 72}$\SP ${238\over 81}$\VL
${77\over 27}$\SP ${1235\over 54}$\SP ${241\over 72}$\SP ${\m907\over 9}$\SP
${\m91\over 6}$\SP ${129\over 1}$\SP ${131\over 8}$\SP ${\m406\over 9}$
\NL $8_8$\VL
${217\over 243}$\SP ${\m1895\over 486}$\SP ${\m283\over 648}$\SP ${589\over 81}$\SP
${25\over 54}$\SP ${\m13\over 3}$\SP ${5\over 24}$\SP ${70\over 81}$\VL
${89\over 81}$\SP ${1621\over 324}$\SP ${193\over 216}$\SP ${\m1012\over 27}$\SP
${\m221\over 36}$\SP ${55\over 1}$\SP ${59\over 8}$\SP ${\m574\over 27}$
\NL $8_9$\VL
${53\over 27}$\SP ${361\over 108}$\SP ${31\over 72}$\SP ${14\over 9}$\SP
${5\over 24}$\SP ${\m28\over 3}$\SP ${\m11\over 12}$\SP ${56\over 9}$\VL
\same
\NL $8_{10}$\VL
${\m707\over 243}$\SP ${\m3593\over 486}$\SP ${\m3007\over 2592}$\SP ${6283\over 324}$\SP
${3229\over 864}$\SP ${\m193\over 9}$\SP ${\m137\over 36}$\SP ${406\over 81}$\VL
${\m10\over 27}$\SP ${2287\over 108}$\SP ${119\over 36}$\SP ${\m1039\over 9}$\SP
${\m53\over 3}$\SP ${155\over 1}$\SP ${79\over 4}$\SP ${\m511\over 9}$
\NL $8_{11}$\VL
${\m238\over 27}$\SP ${\m3049\over 54}$\SP ${\m551\over 72}$\SP ${1613\over 9}$\SP
${313\over 12}$\SP ${\m202\over 1}$\SP ${\m207\over 8}$\SP ${539\over 9}$\VL
${\m2135\over 729}$\SP ${\m1750\over 729}$\SP ${\m3973\over 7776}$\SP ${3793\over 972}$\SP
${5059\over 2592}$\SP ${\m16\over 3}$\SP ${\m61\over 24}$\SP ${\m140\over 243}$
\NL $8_{12}$\VL
${\m4\over 27}$\SP ${661\over 108}$\SP ${253\over 288}$\SP ${\m325\over 36}$\SP
${\m127\over 96}$\SP ${2\over 3}$\SP ${1\over 3}$\SP ${35\over 9}$\VL
\same
\NL $8_{13}$\VL
${\m296\over 243}$\SP ${\m545\over 486}$\SP ${29\over 2592}$\SP ${\m1073\over 324}$\SP
${\m923\over 864}$\SP ${17\over 3}$\SP ${35\over 24}$\SP ${\m119\over 81}$\VL
${\m82\over 81}$\SP ${2521\over 324}$\SP ${1159\over 864}$\SP ${\m5191\over 108}$\SP
${\m2209\over 288}$\SP ${65\over 1}$\SP ${69\over 8}$\SP ${\m637\over 27}$
\NL $8_{14}$\VL
${\m76\over 27}$\SP ${\m2183\over 108}$\SP ${\m719\over 288}$\SP ${1511\over 36}$\SP
${533\over 96}$\SP ${\m34\over 1}$\SP ${\m9\over 2}$\SP ${35\over 9}$\VL
${\m392\over 729}$\SP ${\m4261\over 2916}$\SP ${\m3889\over 7776}$\SP ${5545\over 972}$\SP
${6523\over 2592}$\SP ${\m70\over 9}$\SP ${\m29\over 9}$\SP ${\m35\over 243}$
\NL $8_{15}$\VL
${53\over 3}$\SP ${1331\over 6}$\SP ${1057\over 32}$\SP ${\m3605\over 4}$\SP
${\m4353\over 32}$\SP ${1109\over 1}$\SP ${279\over 2}$\SP ${\m378\over 1}$\VL
${18166\over 2187}$\SP ${72113\over 8748}$\SP ${13781\over 5832}$\SP ${\m17984\over 729}$\SP
${\m21011\over 1944}$\SP ${2449\over 81}$\SP ${1106\over 81}$\SP ${2485\over 729}$
\NL $8_{16}$\VL
${\m5\over 27}$\SP ${377\over 27}$\SP ${73\over 36}$\SP ${\m542\over 9}$\SP
${\m53\over 6}$\SP ${74\over 1}$\SP ${37\over 4}$\SP ${\m224\over 9}$\VL
${\m239\over 81}$\SP ${\m749\over 162}$\SP ${\m169\over 216}$\SP ${259\over 27}$\SP
${181\over 72}$\SP ${\m94\over 9}$\SP ${\m101\over 36}$\SP ${28\over 27}$
\NL $8_{17}$\VL
${\m79\over 81}$\SP ${\m91\over 162}$\SP ${1\over 216}$\SP ${125\over 27}$\SP
${29\over 72}$\SP ${\m28\over 3}$\SP ${\m2\over 3}$\SP ${140\over 27}$\VL
\same
\NL $8_{18}$\VL
${\m146\over 81}$\SP ${\m2347\over 324}$\SP ${\m751\over 864}$\SP ${1975\over 108}$\SP
${613\over 288}$\SP ${\m58\over 3}$\SP ${\m5\over 3}$\SP ${175\over 27}$\VL
\same
\NL $8_{19}$\VL
${\m1486\over 729}$\SP ${\m10007\over 2916}$\SP ${\m547\over 972}$\SP ${2327\over 243}$\SP
${643\over 324}$\SP ${\m329\over 27}$\SP ${\m58\over 27}$\SP ${875\over 243}$\VL
${85\over 3}$\SP ${668\over 3}$\SP ${965\over 32}$\SP ${\m2881\over 4}$\SP
${\m3333\over 32}$\SP ${809\over 1}$\SP ${207\over 2}$\SP ${\m238\over 1}$
\NL $8_{20}$\VL
${22\over 27}$\SP ${779\over 108}$\SP ${23\over 18}$\SP ${\m497\over 9}$\SP
${\m53\over 6}$\SP ${81\over 1}$\SP ${21\over 2}$\SP ${\m287\over 9}$\VL
${253\over 243}$\SP ${\m673\over 243}$\SP ${\m979\over 2592}$\SP ${3175\over 324}$\SP
${1057\over 864}$\SP ${\m11\over 1}$\SP ${\m1\over 1}$\SP ${322\over 81}$
\NL $8_{21}$\VL
${\m52\over 9}$\SP ${\m428\over 9}$\SP ${\m19\over 3}$\SP ${416\over 3}$\SP
${79\over 4}$\SP ${\m147\over 1}$\SP ${\m75\over 4}$\SP ${119\over 3}$\VL
${\m71\over 729}$\SP ${\m6985\over 2916}$\SP ${\m4057\over 7776}$\SP ${7873\over 972}$\SP
${5755\over 2592}$\SP ${\m29\over 3}$\SP ${\m31\over 12}$\SP ${322\over 243}$
&\cr\tabA
}\hrule}
\box1
\endinsert

\midinsert
\noindent Table 1 (continued-1)
\smallskip
\setbox1=\vbox{\offinterlineskip
\hrule
\def\o{&\omit&}
\def\tabA{height3pt
& \omit &
& \omit \o \omit \o \omit \o \omit \o \omit \o \omit \o \omit \o \omit &
& \omit \o \omit \o \omit \o \omit \o \omit \o \omit \o \omit \o \omit & \cr}
\def\NL{&\cr\tabA&\,}
\def\SP{&\omit&}
\def\VL{&&\,}
\def\m{\text{-}}
\def\same{\o\o\o\o\o\o\o}
\halign{&\vrule#&\strut#\hfill\cr
\tabA
& \omit &
&\ $1$\o $\alpha^3$ \o $\alpha^6$ \o $\alpha^2u$ \o $\alpha^5u$ \o $\alpha u^2$ \o $\!\alpha^4\hskip-0.7pt u^2$ \o $u^3$ &
&\ $1$\o $\alpha^3$ \o $\alpha^6$ \o $\alpha^2u$ \o $\alpha^5u$ \o $\alpha u^2$ \o $\!\alpha^4\hskip-0.7pt u^2$ \o $u^3$ &\cr
\tabA
\noalign{\hrule}
\tabA
&\, $9_1$\VL
${\m265\over 3}$\SP $\!\!{\m10993\over 6}$\SP ${\m9047\over 32}$\SP ${28819\over 4}$\SP
${34923\over 32}$\SP ${\m8614\over 1}$\SP ${\m8433\over 8}$\SP $\!{2940\over 1}$\VL
${230278\over 19683}$\SP ${510649\over 39366}$\SP ${23542\over 6561}$\SP $\!\!\!\!{\m260498\over 6561}$\SP
$\!\!\!{\m285377\over 17496}$\SP $\!\!\!{11800\over 243}$\SP $\!\!{39745\over 1944}$\SP $\!\!\!{22519\over 6561}$%
\NL $9_2$\VL
${448\over 9}$\SP ${5303\over 36}$\SP ${157\over 12}$\SP ${\m71\over 3}$\SP
${79\over 8}$\SP ${\m210\over 1}$\SP ${\m177\over 8}$\SP ${553\over 3}$\VL
${\m4463\over 2187}$\SP ${\m45685\over 8748}$\SP ${\m30283\over 23328}$\SP ${54187\over 2916}$\SP
$\!\!{45445\over 7776}$\SP $\!\!{\m1864\over 81}$\SP $\!\!{\m4537\over 648}$\SP $\!\!{1708\over 729}$%
\NL $\bar 9_3$\VL
${134\over 3}$\SP ${\m7367\over 12}$\SP ${\m3599\over 32}$\SP ${14167\over 4}$\SP
${18003\over 32}$\SP ${\m4741\over 1}$\SP ${\m4617\over 8}$\SP $\!{1841\over 1}$\VL
${\m76616\over 6561}$\SP $\!{\m211865\over 13122}$\SP $\!{\m281941\over 69984}$\SP $\!{448729\over 8748}$\SP
$\!\!{415435\over 23328}$\SP $\!\!\!{\m15509\over 243}$\SP $\!{\m42239\over 1944}$\SP ${7315\over 2187}$%
\NL $9_4$\VL
${232\over 3}$\SP ${689\over 12}$\SP ${\m35\over 4}$\SP ${827\over 1}$\SP
${1203\over 8}$\SP ${\m1450\over 1}$\SP ${\m1377\over 8}$\SP ${707\over 1}$\VL
${32555\over 6561}$\SP ${161341\over 26244}$\SP ${140143\over 69984}$\SP ${\m194719\over 8748}$\SP
${\m226657\over 23328}$\SP ${736\over 27}$\SP $\!{2653\over 216}$\SP $\!{6692\over 2187}$%
\NL $\bar 9_5$\VL
${637\over 9}$\SP ${2606\over 9}$\SP ${191\over 6}$\SP ${\m1316\over 3}$\SP
${\m49\over 1}$\SP ${232\over 1}$\SP ${69\over 2}$\SP ${196\over 3}$\VL
${\m6533\over 2187}$\SP ${\m14083\over 2187}$\SP ${\m9679\over 5832}$\SP ${16951\over 729}$\SP
${14725\over 1944}$\SP ${\m2356\over 81}$\SP ${\m740\over 81}$\SP $\!{1960\over 729}$%
\NL $9_6$\VL
${31\over 3}$\SP ${\m2465\over 6}$\SP ${\m2263\over 32}$\SP ${8147\over 4}$\SP
${10203\over 32}$\SP ${\m2618\over 1}$\SP ${\m2541\over 8}$\SP ${980\over 1}$\VL
${\m105118\over 6561}$\SP ${\m289351\over 13122}$\SP ${\m24527\over 4374}$\SP $\!\!{154610\over 2187}$\SP
$\!\!\!{145757\over 5832}$\SP $\!\!\!{\m21332\over 243}$\SP $\!{\m59483\over 1944}$\SP ${7637\over 2187}$%
\NL $9_7$\VL
${94\over 3}$\SP ${\m166\over 3}$\SP ${\m535\over 32}$\SP ${2795\over 4}$\SP
${3759\over 32}$\SP ${\m1046\over 1}$\SP ${\m501\over 4}$\SP ${455\over 1}$\VL
${49886\over 6561}$\SP ${62557\over 6561}$\SP ${206071\over 69984}$\SP ${\m292507\over 8748}$\SP
${\m327541\over 23328}$\SP ${3334\over 81}$\SP ${5717\over 324}$\SP $\!{7049\over 2187}$%
\NL $9_8$\VL
${172\over 243}$\SP ${1003\over 486}$\SP ${209\over 2592}$\SP ${43\over 324}$\SP
${913\over 864}$\SP ${\m43\over 9}$\SP ${\m133\over 72}$\SP ${133\over 81}$\VL
${298\over 81}$\SP ${9575\over 324}$\SP ${3647\over 864}$\SP ${\m12359\over 108}$\SP
${\m4865\over 288}$\SP ${139\over 1}$\SP ${139\over 8}$\SP $\!\!{\m1253\over 27}$%
\NL $9_9$\VL
${124\over 3}$\SP ${\m4849\over 12}$\SP ${\m2443\over 32}$\SP ${9771\over 4}$\SP
${12483\over 32}$\SP ${\m3296\over 1}$\SP ${\m1599\over 4}$\SP $\!{1295\over 1}$\VL
${\m100552\over 6561}$\SP ${\m610997\over 26244}$\SP ${\m432305\over 69984}$\SP ${691145\over 8748}$\SP
$\!{656723\over 23328}$\SP $\!\!{\m23828\over 243}$\SP ${\m33541\over 972}$\SP ${7469\over 2187}$%
\NL $\bar 9_{10}$\VL
${98\over 1}$\SP ${263\over 2}$\SP ${\m3\over 4}$\SP ${732\over 1}$\SP ${1125\over 8}$\SP
${\m1418\over 1}$\SP ${\m1335\over 8}$\SP ${735\over 1}$\VL
${65669\over 6561}$\SP ${70606\over 6561}$\SP ${209935\over 69984}$\SP ${\m287731\over 8748}$\SP
${\m318841\over 23328}$\SP ${1088\over 27}$\SP ${3707\over 216}$\SP $\!{6776\over 2187}$%
\NL $\bar 9_{11}$\VL
${313\over 9}$\SP ${3637\over 18}$\SP ${2471\over 96}$\SP ${\m6467\over 12}$\SP
${\m2393\over 32}$\SP ${546\over 1}$\SP ${567\over 8}$\SP ${\m392\over 3}$\VL
${11050\over 2187}$\SP ${20971\over 4374}$\SP ${7685\over 5832}$\SP ${\m10079\over 729}$\SP
${\m1438\over 243}$\SP ${448\over 27}$\SP ${1609\over 216}$\SP $\!{1435\over 729}$%
\NL $9_{12}$\VL
${329\over 27}$\SP ${9265\over 108}$\SP ${799\over 72}$\SP ${\m2122\over 9}$\SP
${\m787\over 24}$\SP ${240\over 1}$\SP ${123\over 4}$\SP ${\m532\over 9}$\VL
${\m2825\over 729}$\SP ${\m10549\over 2916}$\SP ${\m212\over 243}$\SP ${2083\over 243}$\SP
${1193\over 324}$\SP ${\m308\over 27}$\SP ${\m505\over 108}$\SP ${\m56\over 243}$%
\NL $\bar 9_{13}$\VL
${67\over 1}$\SP ${499\over 4}$\SP ${39\over 8}$\SP ${326\over 1}$\SP
${555\over 8}$\SP ${\m740\over 1}$\SP ${\m339\over 4}$\SP ${420\over 1}$\VL
${61103\over 6561}$\SP ${314719\over 26244}$\SP ${15613\over 4374}$\SP ${\m90109\over 2187}$\SP
${\m49067\over 2916}$\SP ${4096\over 81}$\SP ${6827\over 324}$\SP $\!{6944\over 2187}$%
\NL $9_{14}$\VL
${\m199\over 243}$\SP ${347\over 486}$\SP ${119\over 324}$\SP ${\m868\over 81}$\SP
${\m559\over 216}$\SP ${133\over 9}$\SP ${223\over 72}$\SP ${\m322\over 81}$\VL
${\m307\over 81}$\SP ${\m8675\over 324}$\SP ${\m815\over 216}$\SP ${2804\over 27}$\SP
${553\over 36}$\SP ${\m129\over 1}$\SP ${\m129\over 8}$\SP $\!{1190\over 27}$%
\NL $9_{15}$\VL
${\m1082\over 729}$\SP ${\m3905\over 1458}$\SP ${\m1675\over 1944}$\SP ${2521\over 243}$\SP
${344\over 81}$\SP ${\m374\over 27}$\SP ${\m1157\over 216}$\SP ${49\over 243}$\VL
${491\over 27}$\SP ${3295\over 27}$\SP ${4681\over 288}$\SP ${\m13429\over 36}$\SP
${\m5119\over 96}$\SP ${408\over 1}$\SP ${417\over 8}$\SP $\!\!\!{\m1036\over 9}$%
\NL $\bar 9_{16}$\VL
${7\over 1}$\SP ${\m201\over 1}$\SP ${\m1107\over 32}$\SP ${3751\over 4}$\SP
${4683\over 32}$\SP ${\m1173\over 1}$\SP ${\m561\over 4}$\SP ${434\over 1}$\VL
${\m43018\over 2187}$\SP ${\m255323\over 8748}$\SP ${\m45233\over 5832}$\SP ${71738\over 729}$\SP
${68693\over 1944}$\SP $\!\!\!{\m29651\over 243}$\SP $\!{\m42163\over 972}$\SP ${2597\over 729}$%
\NL $9_{17}$\VL
${313\over 81}$\SP ${3619\over 162}$\SP ${2543\over 864}$\SP ${\m6395\over 108}$\SP
${\m2321\over 288}$\SP ${58\over 1}$\SP ${55\over 8}$\SP ${\m392\over 27}$\VL
${2\over 3}$\SP ${29\over 6}$\SP ${11\over 24}$\SP ${\m29\over 3}$\SP
${\m1\over 6}$\SP ${56\over 9}$\SP ${\m61\over 72}$\SP ${\m7\over 3}$%
\NL $9_{18}$\VL
${52\over 1}$\SP ${75\over 4}$\SP ${\m279\over 32}$\SP ${2415\over 4}$\SP
${3447\over 32}$\SP ${\m1014\over 1}$\SP ${\m120\over 1}$\SP ${483\over 1}$\VL
${83000\over 6561}$\SP ${371311\over 26244}$\SP ${275863\over 69984}$\SP ${\m385519\over 8748}$\SP
${\m419725\over 23328}$\SP ${4390\over 81}$\SP ${3649\over 162}$\SP ${7133\over 2187}$%
\NL $9_{19}$\VL
${26\over 9}$\SP ${541\over 36}$\SP ${211\over 96}$\SP ${\m595\over 12}$\SP
${\m245\over 32}$\SP ${167\over 3}$\SP ${179\over 24}$\SP ${\m49\over 3}$\VL
${652\over 243}$\SP ${2977\over 486}$\SP ${2249\over 2592}$\SP ${\m1565\over 324}$\SP
${\m911\over 864}$\SP ${\m11\over 3}$\SP ${7\over 24}$\SP ${469\over 81}$%
\NL $9_{20}$\VL
${124\over 9}$\SP ${2153\over 36}$\SP ${671\over 96}$\SP ${\m1487\over 12}$\SP
${\m509\over 32}$\SP ${104\over 1}$\SP ${57\over 4}$\SP ${\m35\over 3}$\VL
${13120\over 2187}$\SP ${52589\over 8748}$\SP ${39173\over 23328}$\SP ${\m53933\over 2916}$\SP
${\m59471\over 7776}$\SP ${68\over 3}$\SP ${115\over 12}$\SP ${1183\over 729}$%
\NL $\bar 9_{21}$\VL
${404\over 27}$\SP $\!\!{12997\over 108}$\SP $\!{4669\over 288}$\SP $\!\!{\m13957\over 36}$\SP
${\m5359\over 96}$\SP ${434\over 1}$\SP ${111\over 2}$\SP $\!\!\!{\m1141\over 9}$\VL
${\m3116\over 729}$\SP ${\m15901\over 2916}$\SP ${\m9553\over 7776}$\SP ${15529\over 972}$\SP
${13483\over 2592}$\SP ${\m554\over 27}$\SP ${\m341\over 54}$\SP ${553\over 243}$%
\NL $9_{22}$\VL
${\m172\over 81}$\SP ${667\over 324}$\SP ${79\over 864}$\SP ${\m439\over 108}$\SP
${227\over 288}$\SP ${\m4\over 9}$\SP ${\m65\over 36}$\SP ${\m7\over 27}$\VL
${52\over 81}$\SP ${6689\over 324}$\SP ${2507\over 864}$\SP ${\m7979\over 108}$\SP
${\m3041\over 288}$\SP ${84\over 1}$\SP ${41\over 4}$\SP ${\m707\over 27}$%
\NL $9_{23}$\VL
${21\over 1}$\SP ${12\over 1}$\SP ${\m99\over 32}$\SP ${791\over 4}$\SP
${1167\over 32}$\SP ${\m336\over 1}$\SP ${\m303\over 8}$\SP ${168\over 1}$\VL
${78434\over 6561}$\SP ${201803\over 13122}$\SP ${39467\over 8748}$\SP ${\m114556\over 2187}$\SP
${\m123355\over 5832}$\SP ${5222\over 81}$\SP ${17129\over 648}$\SP ${7301\over 2187}$%
\NL $9_{24}$\VL
${\m94\over 81}$\SP ${2155\over 324}$\SP ${277\over 216}$\SP ${\m1366\over 27}$\SP
${\m607\over 72}$\SP ${215\over 3}$\SP ${59\over 6}$\SP ${\m721\over 27}$\VL
${\m227\over 243}$\SP ${\m1619\over 486}$\SP ${\m967\over 2592}$\SP ${4675\over 324}$\SP
${1405\over 864}$\SP ${\m61\over 3}$\SP ${\m5\over 3}$\SP ${742\over 81}$%
\NL $9_{25}$\VL
${233\over 27}$\SP ${1078\over 27}$\SP ${173\over 36}$\SP ${\m742\over 9}$\SP
${\m253\over 24}$\SP ${67\over 1}$\SP ${69\over 8}$\SP ${\m70\over 9}$\VL
${\m1591\over 729}$\SP ${\m9443\over 2916}$\SP ${\m1997\over 1944}$\SP ${2690\over 243}$\SP
${1603\over 324}$\SP ${\m389\over 27}$\SP ${\m1361\over 216}$\SP ${\m238\over 243}$%
\NL $9_{26}$\VL
${\m917\over 729}$\SP ${\m3094\over 729}$\SP ${\m911\over 972}$\SP ${2938\over 243}$\SP
${2449\over 648}$\SP ${\m121\over 9}$\SP ${\m319\over 72}$\SP ${70\over 243}$\VL
${\m101\over 27}$\SP ${\m3445\over 108}$\SP ${\m307\over 72}$\SP ${838\over 9}$\SP
${161\over 12}$\SP ${\m99\over 1}$\SP ${\m103\over 8}$\SP ${238\over 9}$%
\NL $9_{27}$\VL
${167\over 81}$\SP ${676\over 81}$\SP ${143\over 108}$\SP ${\m970\over 27}$\SP
${\m427\over 72}$\SP ${137\over 3}$\SP ${155\over 24}$\SP ${\m406\over 27}$\VL
${451\over 243}$\SP ${\m541\over 972}$\SP ${\m1\over 162}$\SP ${715\over 81}$\SP
${145\over 216}$\SP ${\m41\over 3}$\SP ${\m17\over 24}$\SP ${574\over 81}$%
\NL $9_{28}$\VL
${\m188\over 27}$\SP ${\m907\over 27}$\SP ${\m155\over 36}$\SP ${706\over 9}$\SP
${131\over 12}$\SP ${\m73\over 1}$\SP ${\m19\over 2}$\SP ${133\over 9}$\VL
${\m2951\over 729}$\SP ${\m20467\over 2916}$\SP ${\m10141\over 7776}$\SP ${17197\over 972}$\SP
${12271\over 2592}$\SP ${\m181\over 9}$\SP ${\m97\over 18}$\SP ${574\over 243}$%
\NL $9_{29}$\VL
${\m38\over 81}$\SP ${1310\over 81}$\SP ${521\over 216}$\SP ${\m2105\over 27}$\SP
${\m415\over 36}$\SP ${100\over 1}$\SP ${99\over 8}$\SP ${\m959\over 27}$\VL
${\m227\over 81}$\SP ${\m283\over 81}$\SP ${\m625\over 864}$\SP ${1309\over 108}$\SP
${943\over 288}$\SP ${\m154\over 9}$\SP ${\m289\over 72}$\SP ${112\over 27}$%
\NL $9_{30}$\VL
${\m4\over 81}$\SP ${901\over 81}$\SP ${1531\over 864}$\SP ${\m5023\over 108}$\SP
${\m2149\over 288}$\SP ${167\over 3}$\SP ${185\over 24}$\SP ${\m469\over 27}$\VL
${\m62\over 243}$\SP ${2159\over 972}$\SP ${1145\over 2592}$\SP ${\m569\over 324}$\SP
${\m743\over 864}$\SP ${\m11\over 3}$\SP ${13\over 24}$\SP ${385\over 81}$%
\NL $9_{31}$\VL
${\m26\over 27}$\SP ${287\over 108}$\SP ${245\over 288}$\SP ${\m2117\over 36}$\SP
${\m923\over 96}$\SP ${95\over 1}$\SP ${95\over 8}$\SP ${\m371\over 9}$\VL
${\m1208\over 729}$\SP ${\m4432\over 729}$\SP ${\m10057\over 7776}$\SP ${18949\over 972}$\SP
${13735\over 2592}$\SP ${\m203\over 9}$\SP ${\m437\over 72}$\SP ${679\over 243}$%
\NL $9_{32}$\VL
${\m61\over 9}$\SP ${\m1469\over 36}$\SP ${\m67\over 12}$\SP ${401\over 3}$\SP
${79\over 4}$\SP ${\m154\over 1}$\SP ${\m20\over 1}$\SP ${140\over 3}$\VL
${\m2981\over 729}$\SP ${\m12391\over 2916}$\SP ${\m1801\over 1944}$\SP ${1918\over 243}$\SP
${2275\over 648}$\SP ${\m82\over 9}$\SP ${\m79\over 18}$\SP ${\m392\over 243}$%
\NL $\bar 9_{34}$\VL
${\m98\over 243}$\SP ${1061\over 972}$\SP ${31\over 81}$\SP ${\m347\over 81}$\SP
${\m175\over 108}$\SP ${3\over 1}$\SP ${7\over 4}$\SP ${133\over 81}$\VL
${19\over 81}$\SP ${722\over 81}$\SP ${1199\over 864}$\SP ${\m3107\over 108}$\SP
${\m1373\over 288}$\SP ${89\over 3}$\SP ${55\over 12}$\SP ${\m182\over 27}$%
\NL $9_{35}$\VL
${772\over 9}$\SP ${3560\over 9}$\SP ${4361\over 96}$\SP ${\m8597\over 12}$\SP
${\m2795\over 32}$\SP ${506\over 1}$\SP ${279\over 4}$\SP ${7\over 3}$\VL
${\m13832\over 2187}$\SP ${\m18799\over 2187}$\SP ${\m47401\over 23328}$\SP ${76165\over 2916}$\SP
${67963\over 7776}$\SP ${\m2650\over 81}$\SP ${\m3431\over 324}$\SP ${1897\over 729}$%
\NL $9_{36}$\VL
${259\over 9}$\SP ${5969\over 36}$\SP ${247\over 12}$\SP ${\m1205\over 3}$\SP
${\m217\over 4}$\SP ${378\over 1}$\SP ${99\over 2}$\SP ${\m224\over 3}$\VL
${5821\over 2187}$\SP ${33725\over 8748}$\SP ${3811\over 2916}$\SP ${\m11393\over 729}$\SP
${\m6301\over 972}$\SP ${514\over 27}$\SP ${439\over 54}$\SP ${1120\over 729}$%
\NL $9_{37}$\VL
${7\over 9}$\SP ${641\over 36}$\SP ${127\over 48}$\SP ${\m361\over 6}$\SP
${\m147\over 16}$\SP ${197\over 3}$\SP ${209\over 24}$\SP ${\m56\over 3}$\VL
${139\over 243}$\SP ${4327\over 486}$\SP ${1705\over 1296}$\SP ${\m2497\over 162}$\SP
${\m1117\over 432}$\SP ${19\over 3}$\SP ${37\over 24}$\SP ${280\over 81}$%
&\cr\tabA
}\hrule}
\box1
\endinsert


\subhead \sectGBZ. The case $\bold k=\Bbb Z$, $\alpha=\beta=0$
\endsubhead

\proclaim{ Proposition \propGBZ } Let $\bold k=\Bbb Z$ and $\alpha=\beta=0$.
We introduce the monomial order on $\bar A$ by saying that
$u^{a_1} v^{b_1} \bar u^{c_1} \bar v^{d_1} > u^{a_2} v^{b_2} \bar u^{c_2} \bar v^{d_2}$
if and only if one of the following conditions holds:
\roster
\item"$\bullet$" either $d_1>d_2$, or $d_1=d_2$ and $c_1>c_2$,
\item"$\bullet$" $(d_1,c_1)=(d_2,c_2)$ and $a_1+b_1>a_2+b_2$,
\item"$\bullet$" $(d_1,c_1,a_1+b_1)=(d_2,c_2,a_2+b_2)$ and $b_1>b_2$.
\endroster
Let
$\Cal G=\{16,4u^2+4v,\,4v^2+4u,\,4uv-4,\,v^3+uv+u^3-3\}$ and
$$
\bar{\Cal G}=\Cal G\cup\{4\bar u-4v,\,4\bar v-4u,\,u\bar u-1,\,v\bar v-1,\,u^3\bar v-3\bar v+v^2+u,\,
\bar u\bar v+u^2\bar v+\bar u v^2-3\}.
$$
Then $\Cal G$ and $\bar{\Cal G}$ are Gr\"obner bases of $I$ and $\bar I$ respectively.
\endproclaim

\noindent{\bf Remark \remMonomOrder }
The monomial order used in Proposition \propGBZ\ can be defined by saying that
this is the lexicographic order with $\bar v>\bar u>w>v'>u'$
under the change of variables $u=uw'$, $v=vw'$.
\medskip

Thus, in the normal form of an element of $\bar A/\bar I$,
the coefficients of the monomials $1$, $u$, $v$ range in $\{0,\dots,15\}$
(note that $u$ and $v$ do not appear in $P(L)$ by Proposition \propPeriod),
the coefficients of
$$
   u^{n+1},\, u^n v,\, u^{n-1}v^2,\, \bar u^n,\, v\bar u^n,\,
   v^2\bar u^n,\, \bar v^n,\, u\bar v^n,\,u^2\bar v^n,\quad  n\ge 1,   
$$
range in $\{0,1,2,3\}$, and all the other coefficients vanish.
Due to Proposition \propInv(h), this fact implies the following nice property of the normal
form of $P(L)$. Let us define the degree of an element of $\bar A/\bar I$
represented by $f\in\bar A$ as
$\min_{g\in f+\bar I} \max_T \deg T$ where $T$ runs over all monomials of $g$ and
$\deg u^av^b\bar u^c\bar v^d\overset{\text{def}}\to=a+b-c-d$.

\midinsert
\noindent Table 1 (continued-2)
\smallskip
\setbox1=\vbox{\offinterlineskip
\hrule
\def\o{&\omit&}
\def\tabA{height3pt
& \omit &
& \omit \o \omit \o \omit \o \omit \o \omit \o \omit \o \omit \o \omit &
& \omit \o \omit \o \omit \o \omit \o \omit \o \omit \o \omit \o \omit & \cr}
\def\NL{&\cr\tabA&\,}
\def\SP{&\omit&}
\def\VL{&&\,}
\def\m{\text{-}}
\def\same{\o\o\o\o\o\o\o}
\halign{&\vrule#&\strut#\hfill\cr
\tabA
& \omit &
&\ $1$\o $\alpha^3$ \o $\alpha^6$ \o $\alpha^2u$ \o $\alpha^5u$ \o $\alpha u^2$ \o $\!\alpha^4\hskip-0.7pt u^2\!\!\!$
 \o $\;\,u^3$ &
&\ $1$\o $\alpha^3$ \o $\alpha^6$ \o $\alpha^2u$ \o $\alpha^5u$ \o $\alpha u^2$ \o $\!\alpha^4\hskip-0.7pt u^2$ \o $u^3$ &\cr
\tabA
\noalign{\hrule}
\tabA
&\,$9_{38}\,$\VL
${125\over 3}$\SP ${1033\over 12}$\SP ${157\over 32}$\SP ${411\over 4}$\SP
${855\over 32}$\SP ${\m304\over 1}$\SP ${\m261\over 8}$\SP ${196\over 1}$\VL
${111548\over 6561}$\SP ${524689\over 26244}$\SP ${48191\over 8748}$\SP $\!\!{\m137809\over 2187}$\SP
${\m146401\over 5832}$\SP ${6278\over 81}$\SP ${20291\over 648}$\SP ${7385\over 2187}$%
\NL $\bar 9_{39}$\VL
${308\over 27}$\SP ${2011\over 27}$\SP ${2857\over 288}$\SP ${\m8437\over 36}$\SP
${\m3223\over 96}$\SP ${261\over 1}$\SP ${267\over 8}$\SP ${\m679\over 9}$\VL
${\m1882\over 729}$\SP ${\m14795\over 2916}$\SP ${\m10757\over 7776}$\SP ${17957\over 972}$\SP
${16763\over 2592}$\SP ${\m635\over 27}$\SP ${\m1715\over 216}$\SP ${371\over 243}$%
\NL $9_{40}$\VL
${\m103\over 27}$\SP ${\m727\over 54}$\SP ${\m503\over 288}$\SP ${1331\over 36}$\SP
${533\over 96}$\SP ${\m41\over 1}$\SP ${\m23\over 4}$\SP ${98\over 9}$\VL
${\m3302\over 729}$\SP ${\m9667\over 2916}$\SP ${\m1759\over 1944}$\SP ${1336\over 243}$\SP
${2467\over 648}$\SP ${\m65\over 9}$\SP ${\m181\over 36}$\SP ${\m749\over 243}$%
\NL $\bar 9_{41}$\VL
${163\over 81}$\SP ${233\over 324}$\SP ${77\over 216}$\SP ${\m176\over 27}$\SP
${\m167\over 72}$\SP ${94\over 9}$\SP ${55\over 18}$\SP ${\m56\over 27}$\VL
${\m61\over 81}$\SP ${\m5789\over 324}$\SP ${\m265\over 108}$\SP ${1709\over 27}$\SP
${325\over 36}$\SP ${\m74\over 1}$\SP ${\m9\over 1}$\SP ${644\over 27}$%
\NL $9_{42}$\VL
${148\over 81}$\SP ${2165\over 324}$\SP ${755\over 864}$\SP ${\m1475\over 108}$\SP
${\m497\over 288}$\SP ${10\over 1}$\SP ${1\over 1}$\SP ${\m35\over 27}$\VL
\same
\NL $9_{43}$\VL
${\m505\over 729}$\SP ${\m553\over 1458}$\SP ${301\over 1944}$\SP ${\m607\over 243}$\SP
${\m205\over 162}$\SP ${3\over 1}$\SP ${13\over 8}$\SP ${182\over 243}$\VL
${41\over 9}$\SP ${1651\over 36}$\SP ${151\over 24}$\SP ${\m460\over 3}$\SP
${\m89\over 4}$\SP ${173\over 1}$\SP ${177\over 8}$\SP ${\m154\over 3}$%
\NL $9_{44}$\VL
${52\over 27}$\SP ${631\over 54}$\SP ${509\over 288}$\SP ${\m1841\over 36}$\SP
${\m755\over 96}$\SP ${65\over 1}$\SP ${67\over 8}$\SP ${\m203\over 9}$\VL
${418\over 243}$\SP ${2705\over 972}$\SP ${1133\over 2592}$\SP ${\m2069\over 324}$\SP
${\m1091\over 864}$\SP ${17\over 3}$\SP ${29\over 24}$\SP ${\m35\over 81}$%
\NL $9_{45}$\VL
${137\over 9}$\SP ${3409\over 36}$\SP ${149\over 12}$\SP ${\m829\over 3}$\SP
${\m313\over 8}$\SP ${295\over 1}$\SP ${303\over 8}$\SP ${\m238\over 3}$\VL
${\m761\over 729}$\SP ${\m5267\over 1458}$\SP ${\m1717\over 1944}$\SP ${3103\over 243}$\SP
${320\over 81}$\SP ${\m425\over 27}$\SP ${\m1019\over 216}$\SP ${406\over 243}$%
\NL $9_{46}$\VL
${\m16\over 9}$\SP ${\m311\over 9}$\SP ${\m491\over 96}$\SP ${1823\over 12}$\SP
${737\over 32}$\SP ${\m194\over 1}$\SP ${\m99\over 4}$\SP ${203\over 3}$\VL
${340\over 243}$\SP ${446\over 243}$\SP ${923\over 2592}$\SP ${\m2399\over 324}$\SP
${\m1313\over 864}$\SP ${82\over 9}$\SP ${59\over 36}$\SP ${\m203\over 81}$%
\NL $\bar 9_{47}$\VL
${\m80\over 27}$\SP ${\m556\over 27}$\SP ${\m889\over 288}$\SP ${3301\over 36}$\SP
${1363\over 96}$\SP ${\m120\over 1}$\SP ${\m31\over 2}$\SP ${385\over 9}$\VL
${\m620\over 243}$\SP ${\m1355\over 486}$\SP ${\m1105\over 2592}$\SP ${709\over 324}$\SP
${859\over 864}$\SP ${\m4\over 3}$\SP ${\m7\over 6}$\SP ${\m119\over 81}$
\NL $9_{48}$\VL
${\m71\over 27}$\SP ${\m397\over 108}$\SP ${\m101\over 144}$\SP ${191\over 18}$\SP
${133\over 48}$\SP ${\m125\over 9}$\SP ${\m233\over 72}$\SP ${28\over 9}$\VL
${95\over 27}$\SP ${3251\over 54}$\SP ${1283\over 144}$\SP ${\m4739\over 18}$\SP
${\m1907\over 48}$\SP ${333\over 1}$\SP ${339\over 8}$\SP $\!\!{\m1036\over 9}$%
\NL $\bar 9_{49}$\VL
${191\over 3}$\SP ${4015\over 12}$\SP ${41\over 1}$\SP ${\m773\over 1}$\SP
${\m825\over 8}$\SP ${705\over 1}$\SP ${741\over 8}$\SP ${\m126\over 1}$\VL
${12389\over 2187}$\SP ${10621\over 2187}$\SP ${8287\over 5832}$\SP ${\m9835\over 729}$\SP
${\m3151\over 486}$\SP ${49\over 3}$\SP ${199\over 24}$\SP ${2366\over 729}$%
&\cr\tabA
}\hrule}
\box1
\endinsert

\proclaim{ Proposition \propSym } If $f(u,v,\bar u,\bar v)$ is the normal form of
an element of $\bar A/\bar I$ of degree $\le2$, then $f(v,u,\bar v,\bar u)$ is the normal form
of the corresponding element. In particular, if $\bar L$ is the mirror image of a link $L$,
then the normal form of $P_{0,0;\Bbb Z}(\bar L)$ is obtained from the normal form of $P_{0,0;\Bbb Z}(L)$ by
swapping $u$ with $v$ and $\bar u$ with $\bar v$.
\endproclaim

Note that if the degree of an element of $\bar A/\bar I$ is greater than two,
then the swapping of $u$ and $v$ can drastically change
the normal form. For example, $u^3$ is already in its normal form whereas
the normal form of $v^3$ is $3uv+3u^3+3$.

In Table 2 we give the normal forms of the invariant $P_{0,0;\Bbb Z}(K)$ for all knots $K$ up to 10 crossings.
We see in this table that there are many repetitions. In Table 3, for each $n=0,1,\dots,12$,
we give the number of different values (up to exchange of $u$ and $v$) that $P_{0,0,;\Bbb Z}$
takes on knots with $\le n$ crossings.

\smallskip\noindent{\bf Example \examFigEight.}
Let us compute $P_{0,0,\Bbb Z}(K)$ where
$K$ is the figure-eight knot $4_1$.
We represent $K$ by the 3-braid
$X=\bar x y\bar x y$ where $x,\bar x,y,\bar y$ are as in Introduction.
Then
{\def\x{\bar x}\def\y{\bar y}
$$
\xalignat2
t(X)&=t\big(\x(-x -y    -\x\y   -\y\x   -\x yx    -x\y x -xy\x)\big) &&\text{by (\eqRelTwo)}           \\
       &=t(    -1 -\x y -\x^2\y -\x\y\x -\x^2yx   -\y x  -y\x)   &&                                     \\
       &=t(    -1 -u\x  -v\x^2  -v\x^2  -u\x\quad -vx    -u\x)   &&\text{by the Markov relation}         \\
       &=t(    -1 -3u\x - 3vx )                                  &&\text{since $\x^2=x$ by (\eqRelOne)} \\
       &=\;\;\,-1 -3uv - 3uv = -1-6uv                            &&\text{by the Markov relation} 
\endxalignat
$$
whence
$$
\xalignat2
 P(K) &= u^{(1-3-0)/2}v^{(1-3+0)/2}t(X)      &&\text{by the definition of $P$}\\
      &= -\u\v - 6                           &&\text{see $t(X)$ computed above}        \\
      &= 7 + u^2\v + \u v^2                  &&\text{Gr\"obner reduction $-\u\v\to u^2\v+\u v^2-3$}
\endxalignat
$$}   

\midinsert
\noindent Table 2. $P_{0,0,\Bbb Z}(K)$ for knots $K$ up to $10$ crossings
\smallskip
\setbox1=\vbox{\offinterlineskip
\hrule
\def\BOL{&}\def\d{}\def\dd{&&}\def\EOL{\ &\cr\noalign{\hrule}\tabAA}\def\b{\overline}
\def\tabA{height1pt & \omit && \omit &\cr}\def\de{\ && \omit &\cr&\ }
\def\tabAA{height1.8pt & \omit && \omit &\cr}\def\de{\ && \omit &\cr&\ }
\halign{&\vrule#&\strut#\hfill\cr
\tabA
&\ $K$  &&\ $P_{0,0,\Bbb Z}(K)$ &\cr
\tabA
\noalign{\hrule}
\tabA
\BOL\ $0_1$\d$6_{2}$\d$6_{3}$\d$8_{6}$\d$8_{7}$\d$8_{8}$\d$8_{9}$\d$8_{16}$\d
$8_{17}$\d$9_{26}$\d$9_{27}$\d$9_{32}$\d$9_{33}$\d$10_{12}$\d$10_{16}$\d$10_{20}$\de
$10_{22}$\d$10_{23}$\d$10_{26}$\d$10_{27}$\d$10_{34}$\d$10_{41}$\d$10_{43}$\d
$10_{48}$\d$10_{52}$\d$10_{54}$\d$10_{79}$\de
$10_{81}$\d$10_{83}$\d$10_{86}$\d$10_{91}$\d$10_{94}$\d
$10_{102}$\d$10_{109}$\d$10_{110}$\d$10_{116}$\de
$10_{118}$\d$10_{123}$\d$10_{125}$\d
$10_{129}$\d$10_{135}$\d$10_{153}$\d$10_{155}$\d$10_{156}$\d$10_{162}$\dd\ $1$\EOL
\BOL\ $\b{3}_{1}$\d$8_{5}$\d$8_{10}$\d$\b{8}_{11}$\d$\b{8}_{20}$\d
$\b{8}_{21}$\d$\b{9}_{24}$\d$\b{9}_{28}$\d$\b{9}_{29}$\d$10_{5}$\d$10_{9}$\d$\b{10}_{32}$\d
$10_{40}$\d$10_{59}$\de
$10_{62}$\d$10_{64}$\d$10_{76}$\d$10_{77}$\d$\b{10}_{82}$\d$10_{84}$\d$\b{10}_{85}$\d
$10_{87}$\d$\b{10}_{103}$\d$10_{106}$\de
$\b{10}_{112}$\d
$10_{113}$\d$\b{10}_{114}$\d$10_{122}$\d
$10_{136}$\d$\b{10}_{141}$\d$\b{10}_{143}$\d$10_{147}$\d$\b{10}_{159}$\dd\ $ \bar u^2 v$\EOL
\BOL\ $4_{1}$\d$8_{4}$\d$8_{13}$\d$9_{22}$\d$9_{30}$\d
$9_{42}$\d$9_{44}$\d$10_{11}$\d$10_{15}$\d$10_{17}$\de
$10_{28}$\d$10_{37}$\d$10_{70}$\d$10_{71}$\d$10_{90}$\d
$10_{93}$\d$10_{104}$\d$10_{119}$\dd\ $7 +  \bar u v^2 + u^2  \bar v$\EOL
\BOL\ $\b{5}_{1}$\d$\b{5}_{2}$\d$\b{7}_{1}$\d$\b{7}_{5}$\d$\b{7}_{6}$\d$\b{8}_{2}$\d$\b{8}_{14}$\d
$9_{3}$\d$\b{9}_{7}$\d$\b{9}_{8}$\d$\b{9}_{9}$\d$\b{9}_{18}$\d$\b{9}_{20}$\d$\b{9}_{31}$\d
$\b{10}_{7}$\d$\b{10}_{8}$\d$\b{10}_{18}$\de
$\b{10}_{24}$\d$\b{10}_{25}$\d$\b{10}_{44}$\d$10_{46}$\d$10_{47}$\d$10_{50}$\d
$10_{51}$\d$10_{56}$\d$10_{57}$\d$10_{92}$\d$10_{95}$\de
$\b{10}_{100}$\d$10_{105}$\d$10_{111}$\d
$10_{117}$\d$\b{10}_{121}$\d$\b{10}_{126}$\d$\b{10}_{127}$\d
$\b{10}_{130}$\d$\b{10}_{131}$\de
$\b{10}_{132}$\d$\b{10}_{133}$\d$\b{10}_{148}$\d$\b{10}_{149}$\d$10_{150}$\d
$10_{151}$\d$\b{10}_{161}$\dd\ $3 + 3  \bar u^3 + 3  \bar u^2 v$\EOL
\BOL\ $\b{6}_{1}$\d$7_{7}$\d$\b{9}_{17}$\d$\b{9}_{34}$\d$10_{4}$\d
$10_{10}$\d$\b{10}_{19}$\d$\b{10}_{29}$\d$\b{10}_{31}$\d$\b{10}_{42}$\de
$\b{10}_{68}$\d$\b{10}_{107}$\d$10_{108}$\d$\b{10}_{146}$\d$10_{158}$\d
$\b{10}_{164}$\dd\ $4 + 3  \bar u^2 v +  \bar u v^2 + u^2  \bar v$\EOL
\BOL\ $\b{7}_{2}$\d$7_{3}$\d$\b{9}_{12}$\d$9_{13}$\d$\b{9}_{25}$\d
$9_{36}$\d$9_{43}$\d$\b{9}_{45}$\d$9_{49}$\d$\b{10}_{2}$\de
$\b{10}_{6}$\d$\b{10}_{38}$\d$\b{10}_{39}$\d$10_{72}$\d$\b{10}_{73}$\d
$10_{157}$\dd\ $2 + 2  \bar u^3 + 2  \bar u^2 v + 3  \bar u^4 v^2$\EOL
\BOL\ $7_{4}$\d$9_{11}$\d$9_{15}$\d$\b{10}_{14}$\d$\b{10}_{21}$\d
$\b{10}_{36}$\d$\b{10}_{67}$\d$10_{69}$\d$\b{10}_{89}$\d
$10_{160}$\d$10_{165}$\dd\ $5 +  \bar u^3 + 3  \bar u^4 v^2$\EOL
\BOL\ $\b{8}_{1}$\d$9_{14}$\d$\b{9}_{41}$\d$\b{10}_{60}$\d$\b{10}_{137}$\d
$10_{138}$\dd\ $13 +  \bar u^3 +  \bar u^2 v +  \bar u v^2 + u^2  \bar v$\EOL
\BOL\ $8_{3}$\d$8_{12}$\d$9_{19}$\d$10_{33}$\d$10_{45}$\d
$10_{88}$\dd\ $15 + 3  \bar u^2 v + 2  \bar u v^2 + 2 u^2  \bar v + 3 u  \bar v^2$\EOL
\BOL\ $\b{8}_{15}$\d$8_{19}$\d$\b{9}_{1}$\d$\b{9}_{6}$\d$9_{16}$\d
$\b{9}_{23}$\d$\b{9}_{38}$\d$\b{10}_{66}$\d$\b{10}_{78}$\d
$10_{139}$\dd\ $3  \bar u^2 v + 3  \bar u^5 v + 3  \bar u^4 v^2$\EOL
\BOL\ $8_{18}$\d$10_{99}$\dd\ $11 + 3  \bar u v^2 + 3 u^2  \bar v$\EOL
\BOL\ $\b{9}_{2}$\d$\b{9}_{4}$\d$9_{10}$\d$\b{10}_{145}$\dd\ $5 +  \bar u^3 + 2  \bar u^2 v +  \bar u^5 v$\EOL
\BOL\ $9_{5}$\dd\ $7 + 3  \bar u^3 + 3  \bar u^2 v +  \bar u^5 v + 3  \bar u^4 v^2$\EOL
\BOL\ $9_{21}$\d$9_{39}$\d$\b{10}_{30}$\dd\ $12 +  \bar u^3 +  \bar u^2 v + 3  \bar u^4 v^2$\EOL
\BOL\ $9_{35}$\dd\ $2 + 2 u  \bar v^2 + 2  \bar v^3 + 2 u^2  \bar v^4 + u  \bar v^5$\EOL
\BOL\ $9_{37}$\dd\ $13 + 3  \bar u^2 v + 3  \bar u v^2 + 3 u^2  \bar v + 3 u  \bar v^2$\EOL
\BOL\ $\b{9}_{40}$\d$10_{61}$\d$10_{65}$\d$\b{10}_{140}$\d$\b{10}_{144}$\d
$10_{163}$\dd\ $12 + 3  \bar u^3 + 2  \bar u^2 v$\EOL
\BOL\ $\b{9}_{46}$\d$9_{47}$\d$10_{75}$\dd\ $15 +  \bar u^3 +  \bar u^2 v$\EOL
\BOL\ $9_{48}$\d$\b{10}_{74}$\dd\ $5 + 2  \bar u^3 + 2  \bar u^2 v$\EOL
\BOL\ $10_{1}$\dd\ $10 +  \bar u v^2 {+\,} u^2  \bar v {\,+\,} 2 u  \bar v^2 {+\,} 2  \bar v^3 {+\,} u^2  \bar v^4$\EOL
\BOL\ $\b{10}_{3}$\d$10_{35}$\dd\ $1 +  \bar u^3 + 2  \bar u v^2 + 2 u^2  \bar v + 3 u  \bar v^2$\EOL
\BOL\ $10_{13}$\d$10_{58}$\dd\ $7 + 3  \bar u^2 v + 3  \bar u v^2 + 3 u^2  \bar v +  \bar v^3$\EOL
\BOL\ $\b{10}_{49}$\d$\b{10}_{53}$\d$\b{10}_{55}$\d$\b{10}_{80}$\d$10_{101}$\d
$10_{124}$\d$10_{128}$\d$10_{134}$\de
$\b{10}_{152}$\d$10_{154}$\dd\ $9 + 2  \bar u^3 {\,+\,}  \bar u^6 {\,+\,} 2  \bar u^2 v {\,+\,} 2  \bar u^5 v {\,+\,}  \bar u^4 v^2$\EOL
\BOL\ $\b{10}_{63}$\d$10_{142}$\dd\ $8 +  \bar u^3 {\,+\,}  \bar u^6 {\,+\,} 2  \bar u^2 v {\,+\,} 3  \bar u^5 v {\,+\,} 2  \bar u^4 v^2$\EOL
\BOL\ $10_{96}$\dd\ $12 +  \bar u^3 + 2  \bar u v^2 + 2 u^2  \bar v$\EOL
\BOL\ $10_{97}$\dd\ $10 + 3  \bar u^3 + 2  \bar u^2 v + 2  \bar u^4 v^2$\EOL
\BOL\ $10_{98}$\dd\ $u^2  \bar v^4$\EOL
\BOL\ $10_{115}$\dd\ $9 + 3  \bar u^2 v +  \bar u v^2 + u^2  \bar v + 3 u  \bar v^2$\EOL
\BOL\ $10_{120}$\dd\ $11 + 3 u  \bar v^2 + u^2  \bar v^4 + u  \bar v^5 +  \bar v^6$ &\cr
}\hrule}
\box1
\endinsert


\midinsert
\noindent Table 3. Here $\Cal K_n$ is the set of knots with $\le n$ crossings; $f(u,v)\sim f(v,u)$.
\smallskip
\vbox{\offinterlineskip
\hrule
\def\o{&\omit&}\def\Card{\operatorname{Card}}\def\K{\Cal K}\def\Z{\Bbb Z}
\def\tabA{height3pt
& \omit &
& \omit \o \omit \o \omit \o \omit \o \omit \o \omit \o
  \omit \o \omit \o \omit \o \omit \o \omit \o \omit \o \omit & \cr}
\def\tabAA{\tabA\noalign{\hrule}\tabA}
\halign{&\vrule#&\strut\hfill\hskip5.25pt#\;\hfill\cr
\tabA
& $n$           && \ 0 \o 1 \o 2 \o 3 \o 4 \o 5 \o 6 \o 7  \o 8  \o 9  \o 10  \o 11  \o 12   &\cr
\tabAA
& $\Card\K_n$   && \ 1 \o 1 \o 1 \o 2 \o 3 \o 5 \o 8 \o 15 \o 36 \o 85 \o 250 \o 802 \o 2978 &\cr
\tabAA
& $\Card P_{0,0;\Z}
 (\K_n)/{\sim}$ && \ 1 \o 1 \o 1 \o 2 \o 3 \o 4 \o 5 \o 7  \o 11 \o 19 \o 29  \o 47  \o 86   &\cr
\tabA
}\hrule}
\endinsert


\subhead\sectHOMFLY. Relation between $P_{0,0;\Bbb Z}$ and the HOMFLY polynomial
\endsubhead
It was discovered by Cabanes and Marin [\refCM] that
$P_{0,0;\Bbb Z/4\Bbb Z}$ is a specialization of the HOMFLY polynomial.
In this subsection we reproduce their arguments in another 
form,
and then we show that if links $L_1$ and $L_2$ have equal HOMFLY polynomials,
than $P_{0,0;\Bbb Z}(L_1)-P_{0,0;\Bbb Z}(L_2)\in\{0,8\}$.
So, $P_{0,0;\Bbb Z}$ can be thought of as a sort of spin of a certain specialization of the HOMFLY polynomial.

Let $k=\Bbb Z/4\Bbb Z$, $K_\infty=K_\infty(0,0;k)$, $A=k[u,v]$, and $I=I(0,0,k)$.
Let $\hat k=k[j]/(j^3-1)$ and let $\hat A=\hat k[u,v]=k[j,u,v]/(j^3-1)$.
Let $H_\infty$ be the Hecke algebra $\hat k B_\infty/(\sigma_1^2+j\sigma_1+j^2)$
(we denote the image of $\sigma_i$ in $H_\infty$ by $g_i$).
A straightforward computation using (\eqRelOne) and (\eqRelTwo)
shows that the correspondence $s_i\mapsto g_i$ defines a ring homomorphism $p:K_\infty\to H_\infty$.
Let $t_H:H_\infty\otimes_{\hat k}\hat A\to M_H\overset{\text{def}}\to=\hat A/(v+ju+j^2)\cong\hat k[u]$
be the Markov trace on $H_\infty$.
The ring $M_H$ is the quotient of $k[j,u,v]$ by the ideal $(j^3-1,v+ju+j^2)$.
For the lexicographic monomial order with $j>w>v'>u'$ ($u=wu'$, $v=wv'$; see Remark \remMonomOrder),
the reduced Gr\"obner base $\Cal G_H$ of this ideal is
$$
           \{v^3+u^3+uv+1,\,   
             j(u^2-v)+uv-1,\,  
             j(uv-1)+v^2-u,\,  
             j(v^2-u)-u^2+v,\, 
             j^2+ju+v\}.       
$$
%
Since $I=(v^3+u^3+uv+1)$, it follows that we can identify $A/I$ with
the subring of $M_H$ generated by $u$ and $v$. Then we have $t=t_H\circ p$, i.~e., the Markov
trace on $K_\infty(0,0;\Bbb Z/4\Bbb Z)$ is determined by the Markov trace on $H_\infty$,
hence $P_{0,0;\Bbb Z/4\Bbb Z}$ is determined by the HOMFLY polynomial.
For example, if we normalize the latter (we denote it $h(a,z)$) by
   $h(\text{unknot})=1$,
   $a h(\,\epsfxsize=8pt\lower1pt\hbox{{\epsfbox{skein-.eps}}}\,)
   -z h(\,\epsfxsize=8pt\lower1pt\hbox{{\epsfbox{skein0.eps}}}\,)
   -a^{-1} h(\,\epsfxsize=8pt\lower1pt\hbox{{\epsfbox{skein1.eps}}}\,)$,
then we have:

\proclaim{ Proposition \propHOMFLY }
$P_{0,0;\Bbb Z/4\Bbb Z}(u,v)=h(ij^2u^{1/2}v^{-1/2},-i)$ where $i=\sqrt{-1}$.
\endproclaim

If we plug $a=ij^2u^{1/2}v^{-1/2}$, $z=-i$ into $h(a,z)$, we obtain a Laurent polynomial in
$i,j,u^{1/2},v^{1/2}$. To get rid of the square roots of $u$, $v$, and $-1$, we just
multiply the result by $(-uv)^{-1/2}$ in the case of a link $L$ with an even number of components
(this corresponds to the normalization (\eqDefP) of $P(L)$).
To eliminate $j$ and to put the result to the canonical form, we reduce it
using the Gr\"obner base $\bar{\Cal G}_H$ of the ideal\break
$(j^3-1,\,v+ju+j^2,\,u\bar u-1,v\bar v-1)$
of the ring $k[j,u,v,\bar u,\bar v]$ for the lexicographic monomial order with
$j>\bar v>\bar u>w>v'>u'$ ($u=wu'$, $v=wv'$; see Remark \remMonomOrder).
We have
$\bar{\Cal G}_H=\Cal G_H\cup
    \{\u u-1,\v v-1,
     u^3\v+\v+v^2+u,\,         
     \u\v+u^2\v+v^2\u+1,\break 
     j(\u-v)-v^2\u+1,\,        
     j(\v-u)+u\v-v \}.         
$
%
%

\smallskip
\noindent{\bf Example \examHOMFLY.} Let $K$ be the trefoil knot  given by the $2$-braid $\sigma_1^{-3}$.
Then $h(K)=-a^4+a^2z^2+2a^2$, hence $P_{0,0;\Bbb Z/4\Bbb Z}(K)=-j^2u^2\v^2-ju\v$
by Proposition \propHOMFLY.
The following code for {\tt Singular} reduces $P(K)$ to $u\bar v^2$.

\smallskip{\tt%
> ring r=(integer,4),(j,V,U,v,u),(lp(3),dp(2));\par
> ideal I=j3-1,v+ju+j2,uU-1,vV-1; reduce(-j2u2V2-juV,std(I));\par
V2u}
\smallskip

\proclaim{ Proposition \propSpin } If $h(L_1)=h(L_2)$, then $P_{0,0;\Bbb Z}(L_1)-P_{0,0;\Bbb Z}(L_2)\in\{0,8\}$.
\endproclaim

\demo{ Proof }
We have shown in Section \sectGBZ\ that the coefficient of all monomials
in the normal form of $P(L)$ except the constant term are in $\{0,1,2,3\}$.
Thus they are determined by $h(L)$ due to Proposition \propHOMFLY.
The constant term is determined mod 8 by the other coefficients due to Proposition \propInv(g).
\qed\enddemo


\Refs
\def\r{\ref}

\r\no\refBF
\by    P.~Bellingeri, L.~Funar
\paper Polynomial invariants of links satisfying cubic skein relations
\jour  Asian J. Math. \vol 8 \yr 2004 \pages 475--509
\endref

\r\no\refCM
\by    M.~Cabanes, I.~Marin
\paper On ternary quotients of cubic Hecke algebras
\jour  Commun. Math. Phys. \vol 314 \yr 2012 \pages 57--92
\endref

\r\no\refFunar
\by    L.~Funar
\paper On cubic Hecke algebras
\jour  Commun. Math. Phys. \vol 173 \yr 1995 \pages 513--558
\endref

\r\no\refMarin
\by    I.~Marin
\paper The cubic Hecke algebra on at most 5 strands
\jour  J. Pure Appl. Algebra \vol 216 \yr 2012 \pages 2754--2782
\endref

\r\no\refO
\by    S.~Yu.~Orevkov
\paper Cubic Hecke Algebras and Invariants of Transversal Links
\jour  Doklady Math. \vol 89 \yr 2014 \pages 115--118
\endref

\r\no\refWebPage
\by    S.~Yu.~Orevkov
\paper Files related to this paper
\jour  http://picard.ups-tlse.fr/~orevkov/fu.html
\endref

\r\no\refOS
\by    S.~Yu.~Orevkov, V.~V.~Shevchishin
\paper Markov theorem for transversal links
\jour  J. of Knot Theory and Ramifications \vol 12 \yr 2003 \pages 905--913
\endref

\endRefs

\enddocument